\documentclass[12pt]{article}

\usepackage{amsfonts}
\usepackage{xr}
\usepackage{color}
\usepackage{graphicx,epsfig}
\usepackage{amsmath}

\def\2{C^{1,2}(\R\times\R^N)}
\def\O{\Omega}
\def\o{\omega}
\def\e{\epsilon}
\def\A{\mathcal{A}}
\def\L{\mathcal{L}}

\def\R{\mathbb{R}}
\def\N{\mathbb{N}}
\def\P{\mathbb{P}}
\def\tilde{\widetilde}
\def\epsilon{\varepsilon}
\def\trait (#1) (#2) (#3){\vrule width #1pt height #2pt depth #3pt}
\def\fin{\hfill\trait (0.1) (5) (0) \trait (5) (0.1) (0) \kern-5pt \trait (5) (5)
(-4.9) \trait (0.1) (5) (0)}

\newtheorem{thm}{\bf Theorem}[section]
\newtheorem{lem}{\bf Lemma}[section]
\newtheorem{prop}[lem]{\bf Proposition}

\newtheorem{defi}[lem]{\bf Definition}

\begin{document}
\title{Spreading speeds for one-dimensional \\
monostable reaction-diffusion equations}

\author{Henri Berestycki$^{\hbox{\small{ a}}}$ and Gr\'egoire Nadin$^{\hbox{\small{
b}}}$\\
\\
\footnotesize{$^{\hbox{a }}$\'Ecole des hautes \'etudes en sciences sociales, CAMS, 190-198 avenue de France, F-75013, Paris, France}\\
\footnotesize{$^{\hbox{b }}$
CNRS, UMR 7598, Laboratoire Jacques-Louis Lions, F-75005 Paris, France}
}

\maketitle


\relpenalty=10000
\binoppenalty=10000
\begin{center}
{\em Dedicated in friendship and admiration to Peter Constantin}
\end{center}
\bigskip
\begin{abstract}
We establish in this article spreading properties for the solutions of equations of
the type $\partial_{t} u -a(x)\partial_{xx}u- q(x)\partial_x u=f(x,u)$,
where $a, q, f$ are only assumed to be uniformly continuous and bounded in $x$, the nonlinearity $f$ is of monostable KPP type 
between two steady states $0$ and $1$ 
and the initial datum is compactly supported.
Using homogenization techniques, we construct two speeds
$\underline{w}\leq \overline{w} $ such that
$\lim_{t\to+\infty}\sup_{0\leq x\leq wt} |u(t,x)-1| = 0$ for all $w\in (0,\underline{w})$ and $\lim_{t\to+\infty} \sup_{x \geq wt} |u(t,x)|
=0$ for all $w>\overline{w}$.
These speeds are characterized in terms of two new notions of
generalized principal eigenvalues for linear elliptic operators in unbounded domains.
In particular, we derive the exact spreading speed when the coefficients are random stationary ergodic,
almost periodic or asymptotically almost periodic 
(where $\overline{w}=\underline{w}$).
\end{abstract}

\noindent {\bf Key-words:} Propagation and spreading properties, Heterogeneous reaction-diffusion equations,
Principal eigenvalues, Linear elliptic operator, Hamilton-Jacobi equations, Homogenization, Random stationarity and ergodicity, Almost
periodicity.

\smallskip

\noindent {\bf AMS classification.} Primary: 35B40, 35B27, 35K57. 
Secondary: 35B50, 35P05, 47B65.

\smallskip

{\em Acknowledgements.} This study was supported in part by the French ``Agence Nationale de la Recherche''
within the project PREFERED.  Henri Berestycki was partially supported by an NSF FRG grant
DMS-1065979.
Part of this work was carried out while the first author was visiting the Department of Mathematics of the University
of Chicago. The authors are grateful to an anonymous referee for useful comments.

\section{Introduction} \label{intro-minspread}

In the present paper, we investigate  the large time behaviour of the solutions of the
Cauchy problem:
\begin{equation} \label{Cauchy} \left\{\begin{array}{rclll}
\partial_{t} u - a(x)\partial_{xx}u-q(x)\partial_x u&=&f(x,u) \quad &\hbox{in}& \quad
(0,\infty)\times\R,\\
u(0,x)&=&u_{0}(x) \quad &\hbox{for all}& \quad x\in\R.\\
\end{array}\right.
\end{equation}
where the coefficients $a$, $q$ and $f$ are assumed to be uniformly continuous and bounded in $x$, with $\inf_{x\in \R} a (x)>0$, the constants $0$ and $1$ are steady states 
of (\ref{Cauchy}) and the reaction term $f$ is supposed
to be monostable between $0$ and $1$. This will be made more precise later in a general framework. 
A typical example of nonlinearity satisfying our hypotheses is
$f(x,s)=c(x)s(1-s)$ with $c$ bounded and $\inf_{\R}c>0$. We consider compactly supported initial 
data $u_0$ with $0\leq u_0\leq 1$ and $u_0\not\equiv 0$. 

We underline that we do not make any other structural assumption in general on the dependence in $x$ of the coefficients, 
such as periodicity, almost periodicity or ergodicity for examples. The results we present is this paper are true for general heterogeneous coefficients. 

We only investigate here the one-dimensional and time-independent case. The multidimensional and space-time heterogeneous 
framework requires a much more involved analysis and will be investigated in our forthcoming 
work \cite{BN2}. 

\bigskip

Equation (\ref{Cauchy}) is indeed a heterogeneous generalization of the classical homogeneous equation
\begin{equation} \label{homogeneous}
\partial_{t}u-\partial_{xx} u = f(u),
\end{equation}
where $f(0)=f(1)=0$ and $f(s)>0$ if $s\in (0,1)$. In the case where $f'(s) \leq f'(0)$ for all $s\in (0,1)$, this equation is called the Fisher-KPP equation 
(from Fisher \cite{Fisher} and Kolmogorov, Petrovsky and Piskunov \cite{KPP}). It is a fundamental equation in  
models of genetics, ecology, chemistry or combustion. 

A central question in these models is to determine precisely
how the steady state $1$ invades the unstable one $0$. 
In the homogeneous framework, Aronson and Weinberger \cite{Aronson} solved this question. They proved that if
$\liminf_{s\rightarrow 0^+} f(s)/s^{1+2/N}>0$, then there exists $w^*>0$ such that
\begin{equation}\label{pbm-spreading} \left\{\begin{array}{lllcl}
\hbox{for all  }w\in (0,w^*),& \ \lim_{t\to+\infty}
\sup_{x\in [0,wt)}& |u(t,x)-1|  &= & 0,\\
\hbox{for all } w> w^*,& \
\lim_{t\to+\infty} \sup_{x\geq wt} &|u(t,x)|  &= & 0.\\
\end{array}\right. \end{equation}
Naturally, a similar result holds for $x\leq 0$. An immediate corollary of (\ref{pbm-spreading}) is that 
$\lim_{t\to +\infty} u(t,x+wt)= 0$ if $w>w^*$ and $\lim_{t\to +\infty} u(t,x+wt)=1$ if $0\leq w<w^*$ locally in $x\in\R$. 
Thus, at large times an
observer moving with speed $w\ge 0$ will only see the stable steady state $1$ if $w < w^*$ and the unstable steady state $0$ if $w > w^*$. We
refer to these results as  {\em spreading  properties} and we call the speed $w^*$ the {\em spreading speed}. It can be proved that $w^*$ is the minimal speed of 
travelling waves
solutions, defined in \cite{Aronson, Fisher, KPP}. If the nonlinear term $f$ is of KPP type, that is, if 
 $f(s)\leq f'(0)s$ for all $s\in (0,1)$, then the spreading speed is explicitely given by $w^*=2\sqrt{f'(0)}$. 
Our goal is to prove spreading properties for the general heterogeneous equation (\ref{Cauchy}). 

Berestycki, Hamel and Nadirashvili \cite{SpeedKPPnote, SpeedKPP2} investigated spreading properties in higher dimension for the homogeneous equation in general unbounded domains such as
spirals, infinite combs or cusps, with Neumann boundary conditions.
In such media, it is not always possible to prove the existence of an exact
spreading speed and several examples are  constructed in \cite{SpeedKPP2} where the spreading speed is infinite or null.
Although our present problem is different from that of \cite{SpeedKPP2}, we expect to observe similar phenomena. As in \cite{SpeedKPP2}, we are
thus led to introduce two speeds:
\begin{equation}\label{def-w^*} \begin{array}{lrlcccl}
w_*:=&\sup \{  w\geq 0,& \lim_{t\to +\infty} \sup_{0 \leq x\leq wt} &|u(t,x)- 1|&=&0 &\},\\
w^*:=&\inf \{ w\geq 0, &\lim_{t\to +\infty} \sup_{x\geq wt}& |u(t,x)|&=&0 &\}.\\
\end{array}\end{equation}
We call these quantities respectively the {\em lower} and {\em upper spreading speeds}. To characterize exactly $w_*$ and $w^*$ in general is still an open problem.
The aim of the present paper is to get accurate estimates on $w_*$ and $w^{*}$. That is, we want to construct $\underline{w}$ (resp. $\overline{w}$) as large as
possible
(resp. as small as possible) such that $\lim_{t\to +\infty} \sup_{0\leq x \leq wt} |u(t,x)-1|=0$ 
for all $w\in (0,\underline{w})$ (resp. $\lim_{t\to +\infty} \sup_{x \geq wt} |u(t,x)|=0$ for all $w>\overline{w}$). In other words, we want to find $\underline{w}$ and $\overline{w}$ as close as possible such that 
$$\underline{w}\leq w_*\leq w^{*}\leq \overline{w}.$$

We underline that the speeds $w_*$ and $w^*$ are not necessarily
equal: there might exist some speeds $w\in (w_*,w^*)$ such that the $\omega-$limit set of $t\mapsto u(t,x+wt)$ is $[0,1]$ for all $x\in\R$. An explicit example has been 
investigated by Garnier, Giletti and the second author \cite{GGN}. 
When $w_*=w^*$, we say that there exists an {\em exact spreading speed}. One of our aim is to identify classes of equations for which there exists an exact spreading speed and to compute it. 

\bigskip

In order to estimate $w_*$ and $w^*$, we will first recall some known spreading properties for some classes of heterogeneous equations.  
A more precise and comprehensive review of known results on this topic will be given in our forthcoming work \cite{BN2}. Here, we just mention the cases of periodic or compactly supported 
spatial heterogeneities, which help to understand our main result. We will also discuss the important class of random stationary ergodic coefficients in Section 
\ref{section-applications} below. 

Consider first a compactly supported heterogeneity: $a\equiv 1$, $q\equiv 0$ and 
$$f(x,s)=\big(b_0-b(x)\big)s(1-s) \quad \hbox{with } b_0>0, \quad b\geq 0, \quad  b \hbox{ smooth and compactly supported.}$$
Then it easily follows from Theorem 1.5 of \cite{BHNa} that $w_*\geq 2\sqrt{b_0}$. It is also easy to check that $w^{*}\leq 2\sqrt{b_0}$ since
$f(x,s)\leq b_0 s (1-s)$. Thus
$$w_*=w^{*}=2\sqrt{b_0}$$
in this case. This example shows that
only what happens at infinity matters to determine $w_*$ and $w^{*}$. 

Next, consider the case where all the coefficients $a$, $q$ and $f$ are periodic in $x$. 
A function $h=h(x)$ is called $L-$periodic if $h(x)=h(x+L)$ for all 
$x\in\R$. The period $L>0$ will be fixed in the sequel. 
It has been proved using various approaches that the spreading property (\ref{pbm-spreading}) still holds in this case \cite{BHNa, Freidlin2, Gartner, Weinberger}. 
When $f$ is of KPP type, that is, when 
$f(x,s)\leq f_s'(x,0)s$ for all $(x,s)\in\R\times [0,1]$,
it is possible to characterize $w^*$ in terms of periodic principal eigenvalues. In this framework, one often speaks of {\em pulled fronts} since 
the propagation speed is determined by the linearization near the unstable steady state $u\equiv 0$. Let now describe the eigenvalues that come 
up in the characterization of the spreading speed. 

Let $\L$ the elliptic operator associated with the linearized equation near $0$:
$$\L\phi:= a(x)\phi''+q(x) \phi'+f'_s(x,0)\phi,$$
and $L_p\phi :=e^{-p x}\L(e^{p x}\phi)$ for all $p\in\R$. Note that $L_p$ has periodic coefficients. Hence, we know from the Krein-Rutman theory that it 
admits a periodic principal eigenvalue $k_p^{per}$, characterized by the existence of a solution $\phi_p$ of
\begin{equation} \label{def-kper}\left\{\begin{array}{rl}
          L_p\phi_p&= k_p^{per}\phi_p,\\
          \phi_p&>0,\\
          \phi_p &\hbox{ is periodic}.\\
         \end{array}\right.\end{equation}
The characterization of the spreading speed $w^*$ proved in \cite{BHNa, Freidlin2, Gartner, Weinberger} then reads
\begin{equation} \label{charw*per} w^*=\min_{p>0} \frac{k_{-p}^{per}}{p}.\end{equation}
This speed is also known to be the minimal speed of pulsating travelling waves \cite{BHNa1}.

\bigskip

Going back to the general case, from these examples we see that in order to find sharp estimates of $w_*$ and $w^{*}$, we need to take into account two aspects.
\begin{itemize}
 \item Only the behaviour of the operator for large $x$ should matter.
\item This behaviour should be characterized through some notion of ``principal eigenvalue'' of the linearized elliptic 
operator near $u=0$. 
\end{itemize}


\section{Statement of the results}\label{section-results}


\subsection{Hypotheses}\label{section-hyp}

We shall assume throughout the paper that $a$, $q$ and 
$f(\cdot,s)$ are uniformly continuous and uniformly bounded with respect to $x\in\R$, uniformly  in $s\in  [0,1]$.
The function $f: \R \times [0,1] \rightarrow \R$ is of class
$\mathcal{C}^{1+\gamma}$ with respect to $s$, uniformly in
$x\in\R$, with $\beta>0$ and $0<\gamma<1$.
We also assume that for all $x\in\R$:
\begin{equation} \label{hyp-f}f(x,0)=f(x,1)=0 \hbox{ and } \inf_{x\in\R} f(x,s)>0 \hbox{ if } s\in
(0,1). \end{equation}
Thus, $0$ and $1$ are steady states of (\ref{Cauchy}). Here, we consider nonlinear terms of KPP type:
\begin{equation} \label{hyp-KPP} f(x,s)\leq f_s'(x,0)s \hbox{ for all }
(x,s)\in\R\times [0,1].\end{equation}
The diffusion coefficient $a$ is supposed to be uniformly positive: $\inf_{x\in\R} a (x)>0$. 
We also require the following condition:
\begin{equation}
\label{hyp-pos}\liminf_{|x|\to +\infty} \Big(4 f_s'(x,0)a(x) -q(x)^2\Big)>0.\end{equation}
This last condition implies, in a sense, that the problem is of monostable (or Fisher-KPP) nature. Indeed,
 if $u_0$ is a non-null initial datum such that $0\leq u_0\leq 1$, 
then by (\ref{hyp-f}) and (\ref{hyp-pos}) the solution $u=u(t,x)$ of (\ref{Cauchy}) converges to $1$ as $t\to +\infty$ 
locally in $x\in\R$ \cite{BHNa, BerestyckiHamelRossi}. In other words, $0$ is an unstable steady state whereas 
the steady state $1$ is globally attractive. 
Note that if $q\equiv 0$ and $f(x,u)= c(x) u (1-u)$, with $\inf_{\R} c (x)>0$, then hypotheses (\ref{hyp-f}), (\ref{hyp-KPP}) and (\ref{hyp-pos}) 
are satisfied. 

These assumptions actually correspond to more general situations with heterogenous steady states $p_-= p_- (x)$ 
and $p_+=p_+ (x)$ instead of $0$ and $1$. Indeed, under the conditions $\inf_\R (p_+-p_-)>0$ and $p_+-p_-$ bounded, the change of variables 
$\tilde{u}(t,x) = \big(u(t,x) - p_- (x)\big) / \big(p_+ (x) -p_-(x)\big)$ reduces the equation with heterogeneous steady states into an equation with steady states $0$ and $1$. 
Thus there is no loss of generality in assuming 
$p_-\equiv 0$ and $p_+ \equiv 1$ as soon as $\inf_\R (p_+-p_-)>0$ and $p_+-p_-$ is bounded.


\subsection{The main tool: generalized principal eigenvalues}\label{section-def}

In order to estimate $w_*$ and $w^*$, we know that one should characterize the heterogeneity of the coefficients 
through some notion of principal eigenvalues associated with 
the elliptic operators defined for all
$\phi\in\mathcal{C}^{2}(\R)$ and $p\in\R$ by
$$\L\phi:= a(x)\phi''+ q(x) \phi'+ f_s'(x,0)\phi \quad \hbox{ and } \quad L_p \phi := e^{-px} \L \big( e^{px} \phi\big).$$
These operators are associated with the linearization near the unstable steady state $0$ of equation (\ref{Cauchy}). As $a$, $q$ and $f$ are just assumed to 
be uniformly continuous
and bounded, with no other assumption such as periodicity for example, these operators are not compact and thus classical eigenvalues do not exist in general. In order to overcome this difficulty,
we need to introduce a generalized notion of principal eigenvalues. 

\begin{defi} \label{defgeneigen}
The {\bf generalized principal eigenvalues} associated with operator $L_p$ in the interval $(R,\infty) \subset\R$, with $R\in \{-\infty\} \cup \R$, are:
\begin{equation}\label{deflambda1'}
\underline{\lambda_1}\big(L_p,(R,\infty)\big):=\sup\big\{\lambda\ |\
\exists\phi\in \A_R \hbox{ such that } L_p\phi\geq\lambda\phi \hbox{ in }(R,\infty) \big\},
\end{equation}
\begin{equation} \label{deflambda1''}\overline{\lambda_1}\big(L_p,(R,\infty)\big):=\inf\big\{\lambda\ |\
\exists\phi\in \A_R \hbox{ such that } L_p\phi\leq\lambda\phi \hbox{ in }(R,\infty) \big\},
 \end{equation}
where if $R\in\R$, $\mathcal{A}_R$ is the set of admissible test-functions over $(R,\infty)$:
\begin{equation} \label{defA} \begin{array}{rll}
 \A_R:= \big\{&\phi\in \mathcal{C}^1\big([R,\infty)\big)\cap\mathcal{C}^{2}\big((R,\infty)\big),&\\ &\phi'/\phi \in L^\infty \big((R,\infty)\big), \
 \phi>0 \hbox{ in } [R,\infty), \ \lim_{x\to +\infty} \displaystyle\frac{1}{x}\ln\phi (x) =0 & \big\}, \\ \end{array}
\end{equation}
and $\A_{-\infty}$ is the set of admissible test-functions over $\R$:
\begin{equation} \label{defAinfty}
 \A_{-\infty}:= \big\{\phi\in \mathcal{C}^{2}(\R),\ \phi'/\phi \in L^\infty (\R), \ \phi>0 \hbox{ in } \R, \ \lim_{|x|\to +\infty} \frac{1}{x}\ln\phi (x) =0 \big\}.
\end{equation}
\end{defi}

Similar quantities have been introduced by Berestycki, Nirenberg and Varadhan \cite{BerNirVaradhan} for multidimensional 
bounded domains with a non-smooth boundary and 
by Berestycki, Hamel and Rossi in \cite{BerestyckiHamelRossi} in unbounded domains (see also \cite{Rossi2, BerRossipreprint}). These 
quantities are involved in the statement of many properties of parabolic and elliptic equations in unbounded domains, such as
maximum principles, existence and uniqueness results \cite{BerestyckiHamelRossi, Rossi2, BerRossipreprint}.
The main difference with 
\cite{BerestyckiHamelRossi, BerNirVaradhan, Rossi2, BerRossipreprint}
is that here we impose $\lim_{x\to +\infty} \frac{1}{x}\ln\phi (x) =0$ instead of asking some bounds from above or below on the test-functions. 
This milder constraint on the test-functions was motivated by the class of random stationary ergodic coefficients, for which 
one can almost surely construct eigenfunctions that are unbounded but satisfy $\lim_{x\to +\infty} \frac{1}{x}\ln\phi (x) =0$.

Although Definition \ref{defgeneigen} is quite simple, these generalized principal eigenvalues are uneasy to handle.  
Several properties of these quantities will be proved in Section \ref{section-geneigen}.
Let only mention here the following result, which proves that if there exists a positive eigenvalue in $\A_R$, 
then the two generalized principal eigenvalues correspond
to the classical notion. This property will be used several times in the sequel to prove spreading properties in periodic, almost periodic and random 
stationary ergodic media. 

\begin{prop} \label{prop-eigen} Let $p\in\R$. 
Assume that there exist $\lambda\in\R$, $R \in \{-\infty\} \cup\R$ and $\phi\in \A_R$ such that $L_p\phi=\lambda\phi$ in $(R,\infty)$. Then,
$$\lambda=\underline{\lambda_1}\big(L_p,(R,\infty)\big)=\overline{\lambda_1}\big(L_p,(R,\infty)\big).$$
\end{prop}


\subsection{Statement of the results}

In order to construct $\underline{w}$ and $\overline{w}$ as precisely as possible such that $\underline{w}\leq w_*\leq w^{*}\leq \overline{w}$, we know that 
only the heterogeneity of the coefficients for large $x$ matters. 
Let thus define for all $p\in\R$:
\begin{equation} \label{defHp-dim1}
\overline{H}(p):=\lim_{R\to +\infty} \overline{\lambda_1}\big(L_p,(R,\infty)\big) \quad\hbox{and}\quad
\underline{H}(p):=\lim_{R\to +\infty} \underline{\lambda_1}\big(L_p,(R,\infty)\big).
\end{equation}
Note that these limits are well-defined since  one can easily prove that $R\mapsto \overline{\lambda_1}\big(L_p,(R,\infty)\big)$ is nonincreasing and 
$R\mapsto \underline{\lambda_1}\big(L_p,(R,\infty)\big)$ is nondecreasing. 
The properties of $\underline{H}$ and $\overline{H}$ are gathered in the next Proposition. 

\begin{prop}\label{prop-Hp}
The functions $\overline{H}$ and $\underline{H}$ are locally Lipschitz-continuous.
Moreover, there exist $C\geq c>0$ such that for all $p \in\R: 
c(1+|p|^2)\leq \underline{H}(p)\leq \overline{H}(p)\leq C(1+|p|^2).$
\end{prop}

We are now in position to define our speeds $\underline{w}$ and $\overline{w}$:
\begin{equation} \label{defw-dim1}
\underline{w}:=\min_{p>0}\frac{\underline{H}(-p)}{p} \quad\hbox{and}\quad
\overline{w}:=\min_{p> 0}\frac{\overline{H}(-p)}{p}.\end{equation}
Our main result is the following.

\begin{thm} \label{spreadingthm}
Take $u_0$ a measurable and compactly supported function such that $0\leq u_0\leq 1$, 
$u_0\not\equiv 0$ and let $u$ the solution of the associated Cauchy problem (\ref{Cauchy}). One has
\begin{enumerate}
\item for all $w>\overline{w}$, $\lim_{t\rightarrow+\infty}\ \sup_{x\geq wt}|u(t,x)|=0,$

\item for all $w\in [0,\underline{w})$, $\lim_{t\rightarrow+\infty}\ \sup_{0\leq x\leq wt}|u(t,x)-1|=0.$
\end{enumerate} 
\end{thm}

In other words, one has $\underline{w}\leq w_*\leq w^{*}\leq \overline{w}$. 
In order to check the optimality of our constructions of $\overline{w}$ and
$\underline{w}$, we will now prove that all the previously known results can be derived from Theorem \ref{spreadingthm}. Moreover, we will show that
 $\underline{w}=\overline{w}$ in various types of media for which no spreading properties have been proved before. 

It is not always true that $\underline{w}= \overline{w}$ since this would imply $w_*=w^*$, that is, the existence of an exact spreading speed. 
But we know from \cite{GGN} that $w_*<w^*$ for some classes of equations. 
Theorem \ref{spreadingthm} would be completely optimal if one was able to prove that $\underline{w}= w_*$ and $\overline{w}= w^*$ are always satisfied, for example by 
proving that $t\mapsto u(t,wt)$ does not converge as $t \to +\infty$ for all $w\in (\underline{w},\overline{w})$. 
We leave this question as an open problem. If Theorem \ref{spreadingthm} was not optimal, then the next step would be 
the improvement of the definitions of our estimates $\underline{w}$ and $\overline{w}$ in order to increase $\underline{w}$ and to decrease $\overline{w}$. 
In particular, we do not know whether our choice of the set of admissible test-functions $\A_R$ is optimal or not. Maybe taking into account more general 
behaviours of the test-functions at infinity could give more accurate estimates on the spreading speeds.


\subsection{Derivation of earlier results}\label{section-applications}

This Section is devoted to some applications of Theorem \ref{spreadingthm}. In particular, we are interested in situations where
$\overline{w}=\underline{w}$. We first prove that when the
heterogeneity is homogeneous, periodic, compactly supported or random stationary ergodic, we recover the known spreading properties.
Then, we show in the next section how to derive new results when the coefficients are almost periodic or asymptotically
almost periodic. We also mention the class of slowly oscillating coefficients, treated in \cite{GGN}. 


\subsubsection*{Homogeneous coefficients}

Assume first that the coefficients are homogeneous, that is, $a$, $q$ and $f$ do not
depend on $x$. Take $a\equiv 1$ and $q\equiv 0$ in order to simplify
the computations. In this case $\L\phi=\phi'' +f'(0)\phi$ and
$L_p\phi=\phi''+2p\phi' +(p^2+f'(0))\phi$.
It immediatly follows from Proposition \ref{prop-eigen} that, taking $\phi\equiv 1$ as a test-function,  
$\underline{\lambda_1}\big(L_p,(R,\infty)\big)=\overline{\lambda_1}\big(L_p,(R,\infty)\big)= p^2+f'(0) \hbox{ for all } R\in\R$. Hence, 
$$\overline{H}(p)=\underline{H}(p)=p^2+f'(0)$$ 
for all $p\in\R$. We conclude that
$\overline{w}=\underline{w}=2\sqrt{f'(0)}.$
This is consistent with Aronson and Weinberger's result \cite{Aronson}.


\subsubsection*{Periodic coefficients}

Assume now that the coefficients are periodic. 
The existence of an exact spreading speed is already known in this case \cite{BHNa, Freidlin2, Gartner, Weinberger}. We now
explain how one can derive this classical result from Theorem \ref{spreadingthm}. 

We know that the operator $L_p$ admits a unique periodic principal eigenvalue
$k_p^{per}$ defined by (\ref{def-kper}), associated with a periodic principal eigenfunction $\phi$. As 
$\phi$ is periodic, continuous and positive, it is bounded and has a positive infimum. Thus $\phi \in \A_R$ for all $R\in\R$
and Proposition \ref{prop-eigen} gives
$\overline{\lambda_1}\big(L_p,(R,\infty)\big) =  \underline{\lambda_1}\big(L_p,(R,\infty)\big)=k_p^{per}$. 
This gives $\overline{H}(p)=\underline{H}(p)=k_p^{per}$ for all $p\in\R$ and
$$\underline{w}=\overline{w}=\min_{p>0}\displaystyle
\frac{k_{-p}^{per}}{p},$$
which is consistent with \cite{BHNa, Freidlin2, Gartner, Weinberger}.


\subsubsection*{Compactly supported heterogeneity}

Assume now that $f_s'(x,0) = b_0 +b(x)$ for all $x\in\R$, where $b_0>0$ and $b$ is a compactly supported and continuous function. Assume that $a\equiv 1$ and 
$q\equiv 0$ in order to simplify the presentation. Assume that $b(x)=0$ for all $|x|\geq r$. Then 
for all $R>r$ and $\phi\in \A_R$, one has $L_p \phi = \phi'' +2p \phi' + (p^2+ b_0)\phi$ in $(R,\infty)$. Thus, Proposition \ref{prop-eigen} gives
$\underline{\lambda_1}\big(L_p,(R,\infty)\big)= \overline{\lambda_1}\big(L_p,(R,\infty)\big)= p^2 +b_0$. Hence $\underline{H}(p) = \overline{H} (p)= p^2 +b_0$
and $$\underline{w} = \overline{w} = w_*=w^*=2\sqrt{b_0}.$$
This is consistent with the result we derived from \cite{BHNa} in the Introduction and even slightly more general since we make no negativity assumption on $b$. 

A generalization of the notion of waves has recently been given by the first author and Hamel \cite{Generalwaves}. One can wonder if the 
speeds $w^*$ and $w_*$ can be viewed as the minimal speed of existence of waves, as in the homogeneous or periodic cases. 
In fact, it has been proved  by Nolen, Roquejoffre, Ryzhik and Zlatos \cite{NolenRoquejoffreRyzhikZlatos} that one can construct some compactly supported heterogeneities
such that the associated equation admits no generalized transition waves. 
Hence, for such heterogeneities, spreading properties hold with speed $w_*=w^*=2\sqrt{b_0}$ but generalized transition waves do not exist.


\subsubsection*{Random stationary ergodic coefficients}

We now consider a probability space $(\O,\P,\mathcal{F})$ and we assume that the reaction rate $f:(x,\omega,s)\in \R\times \O\times [0,1]\to \R$ and the diffusion term $a:(x,\omega) \in \R\times \O \to (0,\infty)$ are random variables.
We assume that $q\equiv a'$, that is, equation (\ref{Cauchy}) is in the divergence form. 
We suppose that $a(\cdot,\omega)$, $a'(\cdot,\omega)$, $1/a(\cdot,\omega)$, $f(\cdot,\omega, s)$ and $f_s'(\cdot,\omega,0)$ are almost surely uniformly
continuous and bounded with respect to $x$ uniformly in $s$, that $f$ is of class $\mathcal{C}^{1+\gamma}$ with respect 
to $s$ uniformly in $x$, that
$f(x,\omega,s) \leq f_s'(x,\omega,0)s$ for all $(x,\omega,s)\in \R\times \O\times [0,1]$. 
The functions $f_s'(\cdot,\cdot,0)$ and $a$ are assumed to be random stationary ergodic. 
This last hypothesis means that there exists a group $(\pi_x)_{x\in\R}$ of measure-preserving
transformations  
such that  $a(x+y,\omega)=a(x,\pi_y \omega)$ and 
$f_s'(x+y,\omega,0)=f_s'(x,\pi_y \omega,0)$ for all $(x,y,\omega)\in \R\times\R\times \O$ and if $\pi_x A=A$ for all $x\in\R$ and for a given  
$A\in\mathcal{F}$, then $\P(A)=0$ or $1$. 
This hypothesis heuristically means that the statistical properties of the medium does not depend on the place where one observes it. 

We expect to compute the speeds 
$\underline{w}$ and $\overline{w}$ for almost every $\omega\in\O$. 
Such a result is already known when the full nonlinearity $f$ (and not only its derivative near $u=0$) is a random stationary ergodic function
since the pioneering work of Freidlin and Gartner \cite{Gartner}. 
They proved that for almost every $\omega\in\O$, one has $w^*=w_*$
and that this exact spreading speed can be computed using a family of Lyapounov exponents associated with the linearization of the equation near $u=0$. 
This result has been generalized by Nolen and Xin to various types of space-time random stationary
ergodic equations \cite{NolenXinrandomgeneral, NolenXinrandom1D, NolenXinrandomshear}. 

Our aim is to check that it is possible to derive $\underline{w}= \overline{w}$ almost surely from Theorem \ref{spreadingthm} and to find a characterization 
of the exact spreading speed that involves the generalized principal eigenvalues. 
The linearized operator now depends on the event $\omega$ and we write for all $\omega\in\O$, $p\in\R$ and $\phi\in\mathcal{C}^2 (\R)$:
\begin{equation}
 L_p^\omega \phi := \big(a(x,\omega)\phi'\big)' +2p a(x,\omega)\phi' + \big(p^2 a(x,\omega) +pa'(x,\omega)+ f_s'(x,\omega)\big)\phi. 
\end{equation}
We associate with these operators two Hamiltonians $\underline{H}^\omega$ , $\overline{H}^\omega$ through (\ref{defHp-dim1}) and two speeds $\underline{w}^\omega$ and 
 $\overline{w}^\omega$ through (\ref{defw-dim1}). 

\begin{prop}\label{rseprop}
Under the hypotheses stated above, there exists a measurable set $\O_0$, with $\P (\O_0)=1$, such that
for all $\omega\in\O_0$:
\begin{equation} \label{charac-rse}\overline{w}^\omega=\underline{w}^\omega=\min_{p>0}\frac{\overline{\lambda_1}(L_{-p}^\omega,\R)}{p}
=\min_{p>0}\frac{\underline{\lambda_1}(L_{-p}^\omega,\R)}{p}\end{equation}
and this quantity does not depend on $\omega\in\O_0$. 
\end{prop}

Hence, the identity $\underline{w}^\omega = \overline{w}^\omega$, which was already known \cite{Freidlin2, Gartner}, can be derived from Theorem \ref{spreadingthm}. Moreover,
we obtain a new characterization of this exact spreading speed involving generalized principal eigenvalues instead the Lyapounov exponents used in \cite{Freidlin2, Gartner}.

The proof of Proposition \ref{rseprop} relies on the equality of the generalized principal eigenvalues.  
\begin{thm} \label{rsethm} There exists a measurable set $\O_0$, with $\P (\O_0)=1$, such that
for all $p\in\R$ and $\omega\in\O_0$:
$$\underline{\lambda_1} (L_p^\omega,\R) = \overline{\lambda_1} (L_p^\omega,\R)$$ 
and this quantity does not depend on $\omega\in\O_0$. 
\end{thm}

The definition of the set of admissible test-functions $\A_{-\infty}$ is important here. If one considers another set of admissible test-functions, for example
$$ \tilde{\A}_{-\infty}= \big\{\phi\in \mathcal{C}^2(\R)\cap L^\infty (\R), \ \inf_{\R}  \phi>0  \big\}$$
then the associated generalized principal eigenvalues are not equal in general. Hence, the class of random stationary ergodic coefficients emphasizes that 
it is very important to use the milder assumption 
$\lim_{|x|\to +\infty} \frac{1}{x}\ln \phi(x)=0$ in the definition of the set of admissible test-functions.


\subsection{Applications to further frameworks and new results}

We have shown in the previous Section how to recover all the previously known results on spreading properties using Theorem \ref{spreadingthm}. 
We will now apply Theorem \ref{spreadingthm} to equations for which no characterization of the exact spreading speed was available up to now. 


\subsubsection*{Almost periodic coefficients}

We will use Bochner's definition of almost periodic functions:

\begin{defi} \cite{Bochner}\label{defalmostper} A function $g:\R\to\R$ is
\emph{almost periodic} with respect to $x\in\R$ if from any sequence
$(x_n)_{n\in\N}$ in $\R$ one can extract a subsequence
$(x_{n_k})_{k\in\N}$ such that $g(x_{n_k}+x)$ converges uniformly
in $x\in\R$.
\end{defi}

\begin{thm}\label{almostperprop}
 Assume that $a$, $q$ and $f_s'(\cdot,0)$ are almost periodic with respect
to $x\in\R$.
Then 
\begin{equation} \label{charac-ap}\overline{w}=\underline{w}=\min_{p>0}\frac{\overline{\lambda_1}(L_{-p},\R)}{p}
=\min_{p>0}\frac{\underline{\lambda_1}(L_{-p},\R)}{p}.\end{equation}
\end{thm}

This proposition is established here as a corollary of the following Theorem about generalized principal eigenvalues, which is new
and of independent interest. 

\begin{thm} \label{almostperthm}
Consider a function $c:\R\to \R$. Assume that $a$, $q$ and $c$ are almost periodic. Let
$\L\phi=a(x)\phi''+q(x)\phi'+c(x)\phi$.
Then one has
$\overline{\lambda_1}(\mathcal{L},\R)=\underline{\lambda_1}(\mathcal{L},\R).$
\end{thm}

It has been shown by Papanicolaou and Varadhan \cite{PapanicolaouVaradhan} that almost periodic functions can be considered as random stationary ergodic ones, 
with an appropriate probability space $(\O,\mathcal{F},\P)$. Roughly speaking, $\O$ is the closure for the uniform convergence of all the translations of 
the almost periodic coefficients (see \cite{PapanicolaouVaradhan} for a precise construction). 
Thus, one could try to apply Theorem \ref{rsethm} in order to derive Theorem \ref{almostperthm}. 
However, one would then get a result for almost every $\omega\in \O$, meaning that this first step may only give the result for a translation at infinity 
of the coefficients. But then, even if such a translation is a good approximation of the original coefficients thanks to the almost periodicity,
as we are investigating large-time behaviours, it is not clear how to control such a behaviour using this approximation. If only 
the reaction term $f$ depends on $x$, then we believe that it is possible to construct appropriate sub and supersolutions and to recover the existence of an exact 
spreading speed using these arguments. But if the diffusion term $a$ and the advection term $q$ depend on $x$, 
it seems that some additional and rather involved arguments should be provided to derive the result not only for a translation at infinity of the operator 
$\mathcal{L}$ but for $\L$ itself. 

We use here a direct approach to prove Theorem \ref{almostperthm}. Namely, we construct appropriate test-functions by using a result 
of Lions and Souganidis \cite{LionsSouganidisap}.


\subsubsection*{Asymptotically almost periodic coefficients}

When the coefficients converge to almost periodic functions at infinity, it is still
possible to prove that $\overline{w}=\underline{w}$. We underline that this result is completely new: the class of asymptotically almost periodic coefficients 
has never been investigated before. 

\begin{prop}\label{asymp-media}
Assume that there exist almost periodic functions $a^*$, $q^*$ and $c^*$ such that
\begin{equation} \label{cv-asympap} \lim_{R\to +\infty}\sup_{x\geq R} \big\{|a(x)-a^*(x)|+|q(x)-q^*(x)|+|f_s'(x,0)-c^*(x)|\big\}=0.\end{equation}
Then $\underline{H}(p)=\overline{H}(p)=\overline{\lambda_1}(L_p^*,\R)=\underline{\lambda_1}(L_p^*,\R)$ for all $p\in\R$ and
\begin{equation}\overline{w}=\underline{w}=\min_{p> 0}\frac{\overline{\lambda_1}(L^*_{-p},\R)}{p}=\min_{p>0}\frac{\underline{\lambda_1}(L^*_{-p},\R)}{p}.
\end{equation}
where $\L^*\phi = a^*(x)\phi'' +q^*(x)\phi' +c^*(x)\phi$ and $L_p^*\phi=e^{-p x}\L^*(e^{p x}\phi)$.
\end{prop}


\subsubsection*{Slowly oscilatting coefficients}

Lastly, let us mention that the method developed in the present paper has been used by Garnier, Giletti and the second author \cite{GGN} to investigate the case where  
$a\equiv 1$, $q\equiv 0$ and $f_s'(x,0) = \mu_0 (\phi (x))$, with $\mu_0$ a periodic function and $\phi$ a smooth increasing function such that $\phi '(x) \to 0$ 
as $x\to +\infty$. If $\phi$ increases sufficiently fast, then these authors proved that 
$w^* = w^{**}$ and it is possible to compute this speed. We refer to \cite{GGN} for a precise definition of ``sufficiently fast''. An example of such a 
$\phi$ is $\phi (x) = \big(\ln (1+|x|)\big)^\alpha, \alpha>1$. This result is proved by constructing appropriate test-functions in the definition 
of the generalized principal eigenvalues (\ref{deflambda1'}). As in the random stationary ergodic setting, it is necessary to consider test-functions which 
are not necessarily bounded but satisfy
$\frac{1}{x}\ln \phi (x)\to 0$ as $|x|\to +\infty$. When $\phi$ increases slowly \big(for example, when 
$\phi (x) =\big(\ln (1+|x|)\big)^\alpha, \alpha \in (0,1)$\big), then it was proved in \cite{GGN} that $w_* = 2\sqrt{\min_\R \mu_0}$ and $w^* = 2\sqrt{\max_\R \mu_0}$, which provides an
example of coefficients for which $w_* < w^*$.


\section{Properties of the generalized principal eigenvalues} \label{section-geneigen}


The aim of this Section is to state some basic properties of the generalized principal eigenvalues and to prove Proposition \ref{prop-Hp}.
In all the Section, we fix an operator $\L$ defined for all
$\phi\in\mathcal{C}^{2}$ by
$$\L\phi:=a(x)\phi''+q(x)\phi'+c(x)\phi,$$
where $a$, $q$ and $c$ are given continuous and uniformly continuous functions over $\R$ and $\inf_{x\in\R} a(x) >0$. We do not require the coefficients to be uniformly continuous in this Section. 

\subsection{Comparison between $\overline{\lambda_1}$ and $\underline{\lambda_1}$}

We begin with an inequality between $\overline{\lambda_1}$ and $\underline{\lambda_1}$.

\begin{prop}\label{prop-comparison}
For all $R\in \{-\infty\} \cup \R$, one has
$$\overline{\lambda_1}\big(\L,(R,\infty)\big)\geq\underline{\lambda_1}\big(\L,(R,\infty)\big).$$
\end{prop}

This comparison result may seem very close to Theorem 1.7 in \cite{BerRossipreprint}. It is not: the test-functions we use here are very different from 
that of \cite{BerRossipreprint}, which were assumed to be either bounded or of positive infimum and to satisfy some boundary condition in $x=R$. 
Here, our condition $\lim_{x\to +\infty} \frac{1}{x} \ln \phi (x)=0$ is milder since $\phi$ could be unbounded and thus the maximum principle of Definition 1.5 
in \cite{BerRossipreprint} does not apply. Moreover, we do not impose any condition at $x=R$ and thus we have no a priori comparison on the test-functions involved in the definitions 
of $\underline{\lambda_1}$ and $\overline{\lambda_1}$. We thus 
need to use a different method to prove the result, which relies on the following technical Lemma. 

\begin{lem}\label{lem:zeps} For all $R\in\R$, there exists no function $z\in \mathcal{C}^2 \big((R,\infty)\big) \cap \mathcal{C}^1 \big([R,\infty)\big)$ 
such that $\lim_{x\to +\infty} \frac{1}{x}\ln z(x) =0$, $z>0$ in $[R,\infty)$ and
\begin{equation}\label{eq:zeps}
 - a(x) z'' - q(x) z'\geq \e z \hbox{ in } (R,\infty), \quad \hbox{with } \e>0.
\end{equation}
\end{lem}

\noindent {\bf Proof.} Assume by translation that $R=0$ in order to enlight the notations. 

\noindent {\bf First case: $z'(0) < 0$.} Let $z_\kappa (x) = e^{\kappa x}z(x)$ 
for all $x>0$, where $\kappa \in \big(0, -z'(0)/z(0)\big)$ will be chosen later. 
This function satisfies
\begin{equation}\label{eq:zkappa}
 - a(x) z_\kappa'' - \big(q(x)-2\kappa a(x)\big) z_\kappa'\geq \big(\e-\kappa q(x)-\kappa^2 a(x)\big) z_\kappa >0 \hbox{ in } (0,\infty),
\end{equation}
taking $\kappa$ small enough. 
Let $m:= \inf_{x\in [0,\infty)} z_\kappa(x)$. As $\lim_{x\to +\infty} z_\kappa(x)=+\infty$ since  $\lim_{x\to +\infty} \frac{1}{x}\ln z(x) =0$, 
this infimum is reached at some $x_0\geq 0$. 
If $x_0>0$,
then a contradiction follows from (\ref{eq:zkappa}). If $x_0=0$, then
$z_\kappa'(0)=z'(0)+\kappa z(0)\geq 0$, which is also a contradiction since our choice of $\kappa$ gives $z'(0)+\kappa z(0)< 0$. 

\smallskip

\noindent {\bf Second case: $z'(0) \geq 0$.} Let $z_\kappa (x) = e^{\kappa x}z(x)$ 
for $\kappa>0$ small enough so that $\e-\kappa q(x)-\kappa^2 a(x)>0$ for all $x\in\R$. 
Equation (\ref{eq:zkappa}) yields that $z_\kappa$ has no local minimum. 
Hence, as $\lim_{x\to +\infty} z_\kappa(x) = +\infty$, the function $z_\kappa$ is nondecreasing. 
Letting $\kappa\to 0$, this gives 
$z'(x)\geq 0$ for all $x>0$. It follows from (\ref{eq:zeps}) that 
\begin{equation}\label{eq:zeps2}
 - z'' +Q z'\geq \delta z \hbox{ in } (0,\infty),
\end{equation}
where $\delta := \e/\sup_{(0,\infty)}a$ and $Q= \max\{ \|q/a\|_\infty,2\sqrt{\delta}+1\}$. 

Let $w_h$ the unique solution of 
$$ -w_h''+Qw_h' = \delta w_h \hbox{ in } (0,h), \quad w_h(0)=z(0), \quad w_h(h) = z(h).$$
This function is explicitly given by 
$$w_h = A_h e^{r_- x} + B_h e^{r_+ x}$$
where $r_\pm = \frac{1}{2} \big( Q \pm \sqrt{Q^2 -4\delta}\big) >0$ and
$$A_h = \displaystyle \frac{ z(h) e^{-r_+ h} -z(0)}{ e^{(r_--r_+)h}-1} \hbox{ and } B_h =z(0)-A_h.$$
As $\lim_{x\to +\infty} \frac{1}{x}\ln z(x) =0$ and $r_+>r_->0$, one has $z(h) e^{-r_+ h} \to 0$ as $h\to +\infty$ and thus $\lim_{h\to +\infty}A_h = z(0)$ and 
$\lim_{h\to +\infty}B_h = 0$. 

Let $\sigma := \inf_{x\in [0,h], w_h(x)\neq 0} z(x)/w_h(x)$.
This quantity is well-defined since $z$ is continuous and positive. Assume by contradiction that $\sigma<1$. Then one has 
$z\geq \sigma w_h$ in $[0,h]$ and there exists $x_0 \in [0,h]$ 
such that $z(x_0)=\sigma w_h (x_0)$. As $z(0)=w_h (0)$, $z(h) = w_h (h)$ and $\sigma<1$, one has $x_0 \in (0,h)$. 
Thus, the function $\zeta := z-\sigma w_h$ is nonnegative and vanishes at the interior point $x_0\in (0,h)$. As it satisfies 
$-\zeta''+Q\zeta'\geq \delta \zeta$, the strong elliptic maximum principle implies $\zeta\equiv 0$. This is a contradiction since 
$\zeta (0) = (1-\sigma) z(0) >0$. Hence $\sigma=1$, which means that $z\geq w_h$ in $(0,h)$.
 
Letting $h\to +\infty$, we eventually get
$$z(x) \geq z(0) e^{r_- x} \hbox{ for all } x>0,$$
which contradicts
$\lim_{x\to +\infty} \frac{1}{x}\ln z(x) =0$.
\hfill $\Box$

\bigskip

\noindent {\bf Proof of Proposition \ref{prop-comparison}.}
Assume by contradiction that $\underline{\lambda_1} \big(\L, (R,\infty)\big) > \overline{\lambda_1} \big(\L, (R,\infty)\big)$, where $R\in\R$. Then there exist $\lambda\in\R$ and $\e>0$ 
such that $$\underline{\lambda_1} \big(\L, (R,\infty)\big) > \lambda>\lambda-\e > \overline{\lambda_1} \big(\L, (R,\infty)\big).$$ This yields that there exist $\phi,\psi \in \A_R$ such that 
$\L \phi \leq (\lambda-\e) \phi$ and $\L \psi \geq \lambda\psi$ in $(R,\infty)$. 
Define $z= \phi/\psi$. Then $z\in \mathcal{C}^2 \big((R,\infty)\big) \cap \mathcal{C}^1 \big([R,\infty)\big)$ is positive over $[R,\infty)$ and satisfies 
$\lim_{x\to +\infty} \frac{1}{x}\ln z(x) =0$. Moreover, one can easily check that $z$ satisfies
$$ a(x) z'' + \Big(q(x)+2 a(x)\frac{\psi' (x)}{\psi (x)} \Big) z' \leq -\e z \hbox{ in } (R,\infty).$$
As $\psi'/\psi \in L^\infty \big((R,\infty)\big)$ since $\psi\in \A_R$, Lemma \ref{lem:zeps} applies and gives a contradiction.

If $R=-\infty$, one has for all $r\in \R$, 
$$\underline{\lambda_1} (\L, \R)\leq \underline{\lambda_1} \big(\L, (r,\infty)\big)\leq \overline{\lambda_1} \big(\L, (r,\infty)\big)\leq  \overline{\lambda_1} (\L, \R).$$
 \hfill $\Box$


\subsection{Continuity with respect to the coefficients}

We now prove that the generalized principal eigenvalues are continuous with respect to the coefficients $q$ and $c$. It is
easy to see that $\underline{\lambda_1}$ and $\overline{\lambda_1}$ are Lipschitz-continuous
with respect
to the zeroth order term of $\L$ using the maximum principle. We improve this property
here and prove that it is also locally Lipschitz-continuous with respect to the
first order term.

\begin{prop} \label{dep-q}
Consider two operators $\L$ and $\L'$ defined for all
$\phi\in\mathcal{C}^{2}$ by
$$\begin{array}{rcl}\L\phi&=&a(x)\phi''+q(x)\phi'+c(x)\phi,\\
\L'\phi&=&a(x)\phi''+r(x)\phi'+d(x)\phi, \\ \end{array}$$
where $a$, $q$, $r$, $c$ and $d$ are continuous and uniformly bounded functions over $\R$ and $\inf_\R a>0$. 
Then, for all $R\in\R$, the following inequalities hold
$$\begin{array}{c} |\overline{\lambda_1}\big(\L',(R,\infty)\big)-\overline{\lambda_1}\big(\L,(R,\infty)\big)|\leq
C\|q-r\|_{L^\infty (R,\infty)} +\|c-d\|_{L^\infty (R,\infty)} + \displaystyle\frac{\|q-r\|_{L^\infty (R,\infty)}^2}{4\inf_{\R} a}, \\
\\
|\underline{\lambda_1}\big(\L',(R,\infty)\big)-\underline{\lambda_1}\big(\L,(R,\infty)\big)|\leq C\|q-r\|_{L^\infty (R,\infty)} +\|c-d\|_{L^\infty (R,\infty)} + 
\displaystyle\frac{\|q-r\|_{L^\infty (R,\infty)}^2}{4\inf_{\R} a} ,\\ \end{array}$$
where $C=\displaystyle\frac{1}{\inf_{\R} a}\max\Big\{\sqrt{\|c\|_{L^\infty (R,\infty)}},\sqrt{\|d\|_{L^\infty (R,\infty)}}\Big\}.$
\end{prop}

\noindent {\bf Remark:} It is an open problem to prove that the generalized principal
eigenvalues are continuous with respect to the diffusion coefficient $a$.

\bigskip

\noindent {\bf Proof.} We use the same type of argument as in the proof of Proposition 5.1 in
\cite{BerNirVaradhan}. Let $\delta=\|q-r\|_{\infty}$ and
$\e=\|c-d\|_{\infty}$. For all constant $M$, one has
$\underline{\lambda_1}(\L+M,(R,\infty))=\underline{\lambda_1}\big(\L,(R,\infty)\big)+M$. Thus, adding a sufficiently large $M$,
one can assume that $c$ and $d$ are
positive functions and that $\underline{\lambda_1}\big(\L,(R,\infty)\big)>0$ and $\underline{\lambda_1}(\L',(R,\infty))>0$.

Take $\kappa>0$. There exists a function $\phi\in \A_R$ and
$$\L\phi=a(x)\phi''+q(x)\phi'+c(x)\phi\geq \Big(\underline{\lambda_1}\big(\L,(R,\infty)\big)-\kappa\Big)\phi
\hbox{ in } (R,\infty).$$
Consider any $\alpha>1$ and define $\psi=\phi^\alpha$. One has $\psi>0$ and $\lim_{x\to +\infty} \frac{1}{x} \ln \psi (x)=0$. Moreover, the function $\psi$ satisfies over $(R,\infty)$:
$$\begin{array}{l}
   -\L'\psi=-a(x)\psi'' -r(x) \psi'-d(x)\psi\\
\\
=-\alpha \phi^{\alpha-1}\big(a(x)\phi''+r(x)\phi'\big)-d(x)\phi^\alpha-\alpha(\alpha-1)\phi^{\alpha-2}
a(x)(\phi')^2 \\
\\
\leq \alpha \delta \phi^{\alpha-1}|\phi'| +\big(\alpha
c(x)-d(x)\big)\phi^\alpha-\Big(\underline{\lambda_1}\big(\L,(R,\infty)\big)-\kappa\Big)\alpha\phi^\alpha -\alpha(\alpha-1)\big(\inf_{ \R} a\big) \phi^{\alpha-2}|\phi'|^2 \\
\\
\leq \displaystyle\frac{\alpha}{4(\alpha-1)\inf_{ \R} a} \delta^2\psi+(\alpha-1) \|c\|_{L^\infty (R,\infty)}\psi+\e
\psi-\Big(\underline{\lambda_1}\big(\L,(R,\infty)\big)-\kappa\Big)\alpha\psi.\\
\end{array}$$
Thus for all $\alpha>1$ and $\kappa>0$, one has:
$$\underline{\lambda_1}\big(\L',(R,\infty)\big)\geq\underline{\lambda_1}\big(\L,(R,\infty)\big)-\kappa-\frac{\alpha}{4(\alpha-1)\inf_{ \R} a}\delta^2-(\alpha-1)
\|c\|_{L^\infty (R,\infty)}-\e.$$
Take $\alpha=1+\delta/\big(2\sqrt{\|c\|_{L^\infty (R,\infty)}\inf_{\R} a}\big)$. Letting
$\kappa\rightarrow 0$, this gives
$$\underline{\lambda_1}\big(\L',(R,\infty)\big)\geq\underline{\lambda_1}\big(\L,(R,\infty)\big)-\delta\sqrt{\frac{\|c\|_{L^\infty (R,\infty)}}{\inf_{ \R} a}}-\e-\frac{\delta^2}{4\inf_{ \R} a}.$$

We get by symmetry:
$$|\underline{\lambda_1}\big(\L',(R,\infty)\big)-\underline{\lambda_1}\big(\L,(R,\infty)\big)|\leq
\delta\max\Big\{\sqrt{\frac{\|c\|_{L^\infty (R,\infty)}}{\inf_{\R} a}},\sqrt{\frac{\|d\|_{L^\infty (R,\infty)}}{\inf_{\R} a}}\Big\}+\e+\frac{\delta^2}{4\inf_{\R} a}.$$
A similar argument, with $0<\alpha<1$, gives the Lipschitz-continuity of
$\overline{\lambda_1}$. \hfill$\Box$


\subsection{Comparison with the classical notion of eigenvalue}

\noindent {\bf Proof of Proposition \ref{prop-eigen}.}
As $\phi \in \A_R$, one can take $\phi$ as a test-function in the definition of $\underline{\lambda_1}$ and $\overline{\lambda_1}$, which gives
$$\lambda \leq \underline{\lambda_1} \big(\L,(R,\infty)\big) \hbox{ and } \lambda\geq \overline{\lambda_1} \big(\L, (R,\infty)\big).$$
On the other hand, Proposition \ref{prop-comparison} yields $\underline{\lambda_1} \big(\L,(R,\infty)\big) \leq \overline{\lambda_1} \big(\L, (R,\infty)\big)$. 
This gives the conclusion. \hfill $\Box$


\subsection{Properties of $\overline{H}$ and $\underline{H}$}

We now gather all the previous results to prove Proposition \ref{prop-Hp}.

\smallskip

\noindent {\bf Proof of Proposition \ref{prop-Hp}.} Recall that, for all $p\in\R$,
\begin{equation} \label{explicitLp}
                   L_p\phi=e^ {-p x}\L (e^ {p x}\phi)=
a(x)\phi'' + \big(2p a(x)+q(x)\big)\phi'+\big(a(x)p^2+q(x)p+f_s'(x,0)\big)\phi.
\end{equation}
Proposition \ref{dep-q} and (\ref{explicitLp}) give the local Lipschitz-continuity
of $\underline{H}$ and $\overline{H}$ with respect to $p$. Proposition \ref{prop-comparison} gives $\overline{H}\geq\underline{H}$.

For all $p\in\R$ and $R>0$, the infimum of the zeroth order term of $L_p$ over $(R,\infty)$ is bounded from below by
$\inf_{x>R} \big(a(x)p^2+q(x) p +f_s'(x,0)\big).$
Thus, taking a constant test-function in the definition of $\underline{\lambda_1}$, one gets
\begin{equation}\label{propHpeq1}\underline{\lambda_1}(L_p,(R,\infty) )\geq \inf_{x>R} \big(a(x)p^2+q(x)p +f_s'(x,0)\big).\end{equation}
Taking the infimum over $p\in\R$ of this inequality, one gets
\begin{equation} \begin{array}{rcl}
\inf_{p\in\R}\underline{\lambda_1}(L_p,(R,\infty) )&\geq& \inf_{ x>R}\inf_{p\in\R} \big(a(x)p^2+q(x)p  +f_s'(x,0)\big)\\
&\geq & \inf_{x>R} \big(f_s'(x,0)-q(x)^2/4 a(x)\big).\\
                 \end{array}
\end{equation}
Eventually, letting $R\to +\infty$ and using (\ref{hyp-pos}), one gets
$$\inf_{p\in\R}\underline{H}(p)\geq \lim_{R\to +\infty}\inf_{x>R} \big(f_s'(x,0)-q(x)^2/4 a(x)\big)>0.$$
Similarly, we obtain from (\ref{propHpeq1}) 
$$\underline{H}(p)\geq \gamma|p|^2 -\|q\|_\infty |p|+ \inf_{x\in\R}f_s'(x,0).$$
Hence, combining these two inequalities, one can find a constant $c>0$ such that for all $p\in\R$,
$$\underline{H}(p)\geq c(1+|p|^2).$$
The other inequality is obtained in a similar way.
\hfill$\Box$


\section{Proof of the spreading property} \label{section-proof}


\subsection{The convergence for $w>\overline{w}$}

We start with the proof of the convergence $\lim_{t\to +\infty} \sup_{x\geq wt} u(t,x)=0$ for all $w>\overline{w}$, which is the easiest part. 

\bigskip

\noindent {\bf Proof of Part 1. of Theorem \ref{spreadingthm}.}
Take $w>\overline{w}$. The definition (\ref{defw-dim1}) of $\overline{w}$ yields that there exist $p>0$ and 
$R>0$ such that $\overline{\lambda_1} \big(L_{-p}, (R,\infty)\big) <wp$. Hence, there exist $\delta \in (0, wp)$ and $\phi\in \A_R$ such that
\begin{equation} \label{eqsupersolR} e^{px} \L ( e^{-px} \phi )  \leq (wp-\delta ) \phi \hbox{ in } (R,\infty).\end{equation}
Up to multiplication of $\phi$ by a positive
constant, one can assume that $\phi (x) e^{-px} \geq u_0 (x)$ for all $x\in\R$ and that $\phi (R)e^{-pR}> 1$ for all $t>0$. 
 Define
$$\overline{u} (t,x) := \left\{\begin{array}{ccc} 1 &\hbox{ if }& x\leq R, \\ 
                                \min \{1, \phi (x) e^{-px + (wp -\delta) t} \} &\hbox{ if }& x>R. \end{array} \right. $$
This function is clearly continuous since $\phi (R) e^{-pR + (wp -\delta)t} \geq \phi (R)e^{-pR} >1$ for all $t\geq 0$.
 
Take now $(t,x)\in (0,\infty)\times \R$ such that  $\overline{u} (t,x) <1$. As $x>R$, it follows from (\ref{eqsupersolR}) that
$$ \partial_t \overline{u} -a(x)\partial_{xx} \overline{u} - q(x)\partial_x \overline{u}-f_s'(x,0)\overline{u} = \Big(-\L ( e^{-px} \phi )+ (wp-\delta) e^{-px} \phi \Big) e^{(wp-\delta)t}\geq 0.$$
Hence, as the constant function $1$ is a supersolution of the Cauchy problem (\ref{Cauchy}), $\overline{u}$ is a weak supersolution of the Cauchy problem. 
The parabolic maximum principle yields $\overline{u} (t,x) \geq u(t,x)$ for all $(t,x) \in (0,\infty) \times \R$. 
It follows from $\lim_{x\to +\infty}\frac{1}{x}\ln \phi (x) = 0$ that 
$$ \sup_{x\geq wt} u(t,x) \leq \sup_{x\geq wt}\big(\phi (x) e^{-px +(wp-\delta) t}\big)
=\sup_{x\geq wt}\big(\phi (x) e^{-px +(wp-\delta)\frac{x}{w}}\big)\leq\sup_{x\geq wt}\big(\phi (x) e^{-\frac{\delta x}{w}}\big) \to 0$$
as $t\to +\infty$. 
\hfill $\Box$


\subsection{The rescaled equation}

In order to prove the convergence for $w\in (0,\underline{w})$ in Theorem \ref{spreadingthm}, we will first determine the limit 
of $v_\epsilon(t,x):=u(t/\e,x/\e)$ as $\e\to 0$ by using homogenization techniques. 
To do so, we follow the ideas developed by Majda and
Souganidis in \cite{MajdaSouganidis}, which are based on 
the half-limits method. 
There is indeed a deep link between homogenization problems and spreading properties for reaction-diffusion equations. This link will be discussed in details
in our forthcoming work \cite{BN2}. We also refer to \cite{LionsSouganidis} for a detailed discussion on homogenization problems and 
on the existence of approximate correctors, a notion which is close to that of generalized principal eigenvalues. 

In the present framework, we need to overcome several new difficulties when trying to apply the ideas of \cite{MajdaSouganidis}
due to the general heterogeneity of the coefficients. 
\begin{itemize}
\item First, classical eigenfunctions of the linearized operators $L_p$ do not exist in general. In
\cite{MajdaSouganidis}, these eigenfunctions play the role of correctors. Here, the
definitions of
$\underline{H}$ and $\overline{H}$ only give sub and super-correctors. As we use
the half-limits method, this difficulty is overcame by using these sub and
super-correctors to prove that
the half-limits are sub and supersolutions of some Hamilton-Jacobi equations. Thus
the generalized principal eigenvalues are well-fitted to our approach.
\item Next, the initial datum becomes $v_\e (0,x) = u_0 (x/\e)$, which depends on $\e$ and converges to $0$ if $x\neq 0$ and to $1$ if $x=0$ when $\e\to 0$. 
This singularity at $x=0$ creates new difficulties since no comparison results exist, as far as we know, for the limit equation (\ref{Z_*subsol}) on the phase $Z_*$. 
We will overcome this difficulty by getting estimates on the phase $Z_*$ by hand. 
\item Lastly, we want to prove that only what happens in $(R,\infty)$ with large $R$ plays a role in the
computations of $\underline{H}$ and $\overline{H}$.
\end{itemize}

The authors are indebted to Jimmy Garnier and Thomas Giletti for a careful reading and useful comments on this proof. 

\bigskip

The proof starts as in \cite{MajdaSouganidis}.
Define
\begin{equation} \label{defZepsilon}
Z_\e (t,x) := \e \ln v_\e (t,x) \end{equation}
and let
\begin{equation} \label{defsubsuper}
Z_*(t,x):=\liminf_{(s,y)\rightarrow (t,x), \e\rightarrow 0} Z_\e(s,y).
\end{equation}
Our aim is to check that $Z_*$ is a supersolution of some Hamilton-Jacobi equation to be determined. 

The following Lemma ensures that the function $Z_*$ takes finite values for all $t>0, x>0$. 

\begin{lem}\label{prop-Z*}
 The family $(Z_\e)_{\e>0}$ satisfies the following properties:

i) For all $t>0$, one has $Z_*(t,0)=0$. 

ii) For all compact set $Q\subset (0,\infty)\times \R$, there exist a constant $C=C(Q)$ and $\e_0=\e_0(Q)$ such that $|Z_\e (t,x)|\leq C$ for all $0<\e<\e_0$ and $(t,x)\in Q$. 
\end{lem}

{\bf Proof.}

i) We know from \cite{BHNa} that, as (\ref{hyp-pos}) is satisfied, there exists $c>0$ such that
\[\lim_{t\rightarrow+\infty}\inf_{|x|\leq ct} u(t,x)=1.\]
Fix $t_0>0$ such that $\inf_{|x|\leq ct} u(t,x)\geq 1/2$ for all $t\geq t_0$.
Consider now $t>0$ and a sequence $(s_n,y_n)\in\R^+\times\R^N$ such that $s_n\rightarrow t$ and $y_n\rightarrow 0$ as $n\rightarrow +\infty$. 
Thus $|y_n|/s_n\leq c$ and $s_n/\e \geq t_0$ when $n$ is large and $\e$ is small. This yields
\[0\geq Z^\e(s_n,y_n)=\e\ln u(s_n/\e,y_n/\e)\geq \e \ln \inf_{|x|\leq c s_n/\e} u(s_n/\e,x)\geq -\e \ln 2 \rightarrow 0 \hbox{ as } \e\rightarrow 0.\]
Thus $Z_*(t,0)=0$.

\smallskip

ii) First, the Krylov-Safonov-Harnack inequality \cite{KrylovSafonov} yields that for any $R\leq 2$, there exists a constant $C>0$, 
such that
$$\forall |y|\leq R/2, \quad \forall s>0, \quad u(s+R^2,0)\leq  C u(s+2 R^2,y).$$
Iterating this inequality, one gets for all $N\in\mathbb{N} \backslash \{0\}$:
\begin{equation} \label{eq:HarnackN} \forall |y|\leq N R/2, \quad u(N R^2,0)\leq  C^N u(2 N R^2,y).\end{equation}

Let now derive the local boundedness of $Z_*$ from this inequality. Take $T>\tau>0$, $\rho>0$, $t\in (\tau,T)$ and $|x|\leq \rho$. 
Let $R_0 := \tau /8 \rho$ and
$N_\e:= [t/(2\e R_0^2)]$ the integer part of $t/(2\e R_0^2)$ for any $\e \in (0,1)$. Define $R_\e:= \sqrt{t/(2\e N_\e)}$. It is easily noticed 
that $\lim_{\e \to 0^+} R_\e = R_0$ using the definition of the integer part.
Let $y= x/\e$. One has 
$$|y| =\frac{|x|}{\e} \leq \frac{\rho}{\e} = \frac{\tau}{8 \e R_0} \leq \frac{\tau }{4\e R_\e}$$
for $\e$ small enough, from which we get $|y|\leq N_\e R_\e /2$ since $t\geq \tau$ and $N_\e R_\e^2 = t/2\e$ by construction. 
Hence, we can apply inequality (\ref{eq:HarnackN}), which gives:
$$ u \big(t/(2\e),0\big)=u(N_\e R_\e^2,0)\leq  C^{N_\e} u(2 N_\e R_\e^2,y)=C^{N_\e} u (t/\e,x/\e).$$
As $Z_\e(t,x)=\e \ln u(t/\e,x/\e)$, we get 
\begin{equation} \label{eq:HarnackZe} Z_\e (t/2,0) \leq \e N_\e \ln C + Z_\e (t,x).\end{equation}
Moreover, we know that 
$$\lim_{\e \to 0^+} \e N_\e =\lim_{\e \to 0^+} \frac{t}{2 R_\e^2} = \frac{t}{2 R_0^2} \leq \frac{T}{2R_0^2} = \frac{32 T \rho^2}{\tau^2}.$$ 
Hence, there exist $\e_0>0$ and a constant, that we still denote $C$, such that for all $\e\in (0,\e_0)$, $t\in [\tau, T]$ and $|x|\leq \rho$, one has 
$$|Z_\e (t,x) | = -Z_\e (t,x) \leq C-Z_\e (t/2,0).$$
Step i) yields that $\lim_{\e \to 0^+} Z_\e (t/2,0)=0$ locally uniformly in $t\in (\tau,T)$,  
which ends the proof. 
\hfill$\Box$

\bigskip

The link between $Z_*$ and the convergence of $v_\e$ as $\e\to 0$ is given by the following Lemma:
\begin{lem} \label{cvHJ}
One has 
\begin{equation} \liminf_{\e\to 0} v_\e (t,x)>0
\hbox{ for all }
(t,x) \in {\rm int} \{Z_*=0\}.
\end{equation}
\end{lem}

\noindent {\bf Proof.}
Fix $(t_0,x_0)\in {\rm int}\{Z_*=0\}$. As $u(t,x)\to 1$ as $t\to +\infty$ locally in $x$, one has $v_\e(t,0)\to 1$ as $\e\to 0$ for all $t>0$. 
We thus assume that $x_0 \neq 0$. One has $Z^\e(t,x)\to 0$ as $\e\to 0$ uniformly in the neighborhood of $(t_0,x_0)$.
Define
\[\phi(t,x)=-|x-x_0|^2-|t-t_0|^2.\]
As $Z_*=0$ in the neighbourhood of $(t_0,x_0)$ and $\phi$ is nonpositive, the function $Z_{\e}-\phi$ reaches a minimum at a point $(t_\e,x_\e)$, with
$(t_\e,x_\e)\to (t_0,x_0)$ as $\e\to 0$.
Thus, the equation on $Z_\e$ (see (\ref{eqZepsilon}) below) gives
\[\partial_t\phi-\e a\partial_{xx}\phi-a(x_\e/\e) (\partial_{x}\phi)^2-q(x_\e/\e)\partial_{x}\phi-(v_{\e})^{-1} f(x_\e/\e,v_{\e})\geq 0,\]
where the derivatives of $\phi$ and $v_\e$ are evaluated at
$(t_\e,x_\e)$. An explicit computation of the left-hand side gives
\[(v_{\e})^{-1} f(x_{\e}/\e,v_{\e}(t_{\e},x_{\e}))\leq o(1) \hbox{ at } x_{\e} \hbox{ as } \e\rightarrow 0.\]
As $f$ is of class $\mathcal{C}^{1+\gamma}$ with respect to $s$ uniformly in $x$, there exists $C>0$ such that for all
$(x,u)\in\R\times [0,1]$,
\[ f(x,u)\geq f_s'(x,0)u -Cu^{1+\gamma}. \]
This gives
\[f_s'(x_{\e}/\e,0) \leq Cv_{\e}(t_{\e},x_{\e})^{\gamma}+o(1) \hbox{ as } \e\to 0.\]
Hypothesis (\ref{hyp-pos}) together with $x_0>0$ give
\[\liminf_{\e\to 0} f_s'(x_{\e}/\e,0)>0.\]
Thus $\liminf_{\e\rightarrow 0} v_\e(t_\e,x_\e)>0$. 

Next, the definition of $(t_\e,x_\e)$ yields
$$Z_{\e}(t_0,x_0)=Z_{\e}(t_0,x_0)-\psi (t_0,x_0) \geq Z_{\e}(t_\e,x_\e)-\psi (t_\e,x_\e)\geq Z_{\e}(t_\e,x_\e).$$
It follows from the definition of $Z_\e$ that $$\e\ln u (t_0/\e,x_0/\e ) =\e \ln v_\e (t_0,x_0)\geq \e\ln u(t_\e/\e,x_\e/\e)=\e \ln v_\e (t_\e,x_\e).$$
Hence, $v_\e (t_0,x_0)\geq v_\e (t_\e,x_\e)$ and one gets $\liminf_{\e\to 0} v_\e (t_0,x_0) >0$. 
\hfill$\Box$


\subsection{The equation on $Z_*$}

In order to identify the set ${\rm int} \{Z_*=0\}$, we prove in this Section that $Z_*$ is a supersolution of some first order Hamilton-Jacobi equation. 
As $u$ satisfies (\ref{Cauchy}), the definition of $Z_\e$ yields for all $\e>0$:
\begin{equation} \label{eqZepsilon}
\partial_{t} Z_\e - \e a(x/\e) \partial_{xx} Z_\e- a(x/\e) \big(\partial_x Z_\e\big)^2- q(x/\e)\partial_x Z_\e=\frac{1}{v_\e}f(x/\e,v_\e).
\end{equation}
The keystone of our proof is the next result.
\begin{prop}\label{eqZ^*Z_*}
The function $Z_*$ is a lower semi-continuous viscosity solution of
 \begin{equation} \label{Z_*subsol} \left\{
\begin{array}{l}
    \max \{ \partial_t Z_*- \underline{H}(\partial_x Z_*), Z_* \}\geq 0 \hbox{ in } (0,\infty)\times (0,\infty),\\
Z_*(t,0)=0 \hbox{ for all } t>0.\\
\end{array}\right. \end{equation}
\end{prop}

\noindent {\bf Proof.} As $Z_*\leq 0$ since $u\leq 1$, we need to prove that
\begin{equation} \label{Z_*equation1}\partial_t Z_* - \underline{H}(\partial_{x} Z_*)\geq 0 \hbox{ in } \{Z_*<0\}.\end{equation}
Fix a smooth test function $\phi$ and assume that $Z_*-\phi$ admits a strict minimum at some point $(t_0,x_0)\in (0,\infty)\times (0,\infty)$ over the ball
$\overline{B_r}:= \{(t,x)\in (0,\infty)\times (0,\infty), |t-t_0|+|x-x_0|\leq r\}$, with $Z_* (t_0,x_0)<0$. 
Define $p=\partial_x\phi(t_0,x_0)$.
If we manage to prove that for all $\mu>0$,
\[\partial_t\phi(t_0,x_0)-\underline{H}(p)\geq -\mu,\]
then letting $\mu\to 0^+$ would imply that $Z_*$ is a viscosity subsolution of (\ref{Z_*equation1}). 

Fix $R>0$ and consider a function $\psi\in\A_R$ such that $\Big(L_p-\underline{\lambda_1}\big(L_p,(R,\infty)\big)\Big)\psi\geq -\mu\psi$. Let
$w=\ln \psi$, this function satisfies over $(R,\infty)$:
\begin{equation} \label{Z^*eqw} -a(x) \big(w''+(w' +p)^2\big)-q(x)\big(w'+p\big)\leq f_s'(x,0)
-\underline{\lambda_1}\big(L_p,(R,\infty)\big)+\mu.\end{equation} 
Moreover, one has $\e w(x/\e) \to 0$ as $\e\to 0$ locally in $x\in (R,\infty)$. 

The definition of $Z_*$ yields that there exist a sequence of positive numbers $(\e_n)_n$ and a sequence $(s_n,y_n)_n$ in $\overline{B_r}$ such that $\e_n\to 0$, 
$s_n\to t_0$, $y_n \to x_0$
and $Z_{\e_n}(s_n,y_n)\to Z_*(t_0,x_0)$ as $n\to +\infty$. 
For all $n$, let $(t_n,x_n)\in \overline{B_r}$ such that the function 
\begin{equation}
 Z_{\e_n}-\phi-\e_n w(\cdot/\e_n) \hbox{ reaches a minimum at } (t_n,x_n) \hbox{ over } \overline{B_r}.
\end{equation}
As the sequence $(t_n,x_n)_{n}$ lies in $\overline{B_r}$ one can assume, up to extraction, that it converges in $\overline{B_r}$. 
Let $(T_0,X_0)$ its limit. 
For all $n$ and for all $(t,x)\in \overline{B_r}$, one has 
\begin{equation} \label{eqZt_n} Z_{\e_n}(t,x)-\phi (t,x)-\e_n w(x/\e_n)\geq Z_{\e_n}(t_n,x_n)-\phi (t_n,x_n)-\e_n w(x_n/\e_n).\end{equation}
Taking $t=s_n$, $x=y_n$ and letting $n\to +\infty$, the definition of $Z_*$ yields that, 
$$ Z_*(t_0,x_0)-\phi (t_0,x_0)\geq Z_*(T_0,X_0)-\phi (T_0,X_0) \hbox{ for all } (t,x)\in \overline{B_r}.$$
Hence, as $Z_*-\phi$ reaches a strict local minimum at $(t_0,x_0)$ over the ball $\overline{B_r}$, one gets $(T_0,X_0)=(t_0,x_0)$. 
We have thus proved that
\begin{equation} \begin{array}{l}
Z_{\e_n}(t_n,x_n)\rightarrow Z_*(t_0,x_0),\\
(t_n,x_n)\rightarrow (t_0,x_0) \hbox{ as } n\rightarrow+\infty,\\
Z_{\e_n}-\phi-\e_n w(\cdot/\e_n) \hbox{ reaches a local minimum at } (t_n,x_n).\\
\end{array} \end{equation}

As $x_0> 0$, one has $x_n/\e_n\rightarrow+\infty$. 
Take $n$ large enough such that $x_n/\e_n>R$. As $Z_{\e_n}-\big(\phi+\e_n w
(\frac{\cdot}{\e_n})\big)$ reaches a
local minimum in $(t_n,x_n)$, we get:
\begin{equation}\label{Z^*estimate0} \begin{array}{l}
\partial_t\phi-\partial_t Z_{\e_n}-\e_n a (x_n/\e_n)(\partial_{xx}\phi+\e_n^{-1}\partial_{xx} w-\partial_{xx} Z_{\e_n})\\
- a(x_n/\e_n)(\partial_x\phi+\partial_{x} w-\partial_{x} Z_{\e_n})^2-q(x_n/\e_n)(\partial_{x}\phi+\partial_{x} w-\partial_{x} Z_{\e_n})\geq 0,\\
\end{array} \end{equation}
where the derivatives of $\phi$ and $Z_{\e_n}$ are evaluated at
$(t_n,x_n)$ and the derivatives of $w$ are evaluated at $x_n/\e_n$.
Using the equation (\ref{eqZepsilon}) satisfied by $Z_\e$, we get
\begin{equation}\label{Z^*estimate00}\begin{array}{l}
\partial_t\phi - a(x_n/\e_n)(\e_n\partial_{xx}\phi+\partial_{xx}w)- a(x_n/\e_n)(\partial_{x}\phi+\partial_{x} w)^2
-q(x_n/\e_n)(\partial_{x}\phi+\partial_{x} w)\\
\geq \frac{1}{v_{\e_n}}f(x_n/\e_n,v_{\e_n}).\\ \end{array}\end{equation}
As $Z_{\e_n}(t_{n},x_{n})\rightarrow Z_*(t_0,x_0)<0$, one has 
\[v_{\e_n}(t_{n},x_{n})= \exp [ \e_n^{-1} Z_{\e_n}(t_{n},x_{n})]\rightarrow 0 \hbox{ as } n\rightarrow+\infty\]
and thus the right-hand side of equation (\ref{Z^*estimate00}) is equivalent to $f_s'(x_n/\e_n,0)$ as $n\to +\infty$. 
Using (\ref{Z^*eqw}), we deduce
\[\begin{array}{rl}
\partial_t\phi-\overline{\lambda_1}\big(L_p,(R,\infty)\big) \geq &-\mu+\e_n a(x_n/\e_n) \partial_{xx}\phi+q(x_n/\e_n) (\partial_{x}\phi-p)\\
&-\|a\|_\infty(\partial_{x}\phi-p)^2- 2\|a\|_\infty |\partial_{x}\phi-p||\partial_{x} w +p|+o(1),\\ \end{array}\]
where the derivatives of $\phi$ are evaluated at
$(t_n,x_n)$. We remind to the reader that $p=\partial_{x}\phi(t_0,x_0)$.
Hence, letting $n\rightarrow +\infty$ and $\mu \to 0$, this leads to
\[\partial_t\phi(t_0,x_0) -\underline{\lambda_1}\big(L_p,(R,\infty)\big)\geq 0.\]

Finally, letting $R\to+\infty$ , one has 
\begin{equation} \label{Z^*equation} \max \{ \partial_t Z_* - \underline{H}(\partial_x Z_*), Z_* \}\geq 0 \hbox{ in } (0,\infty)\times (0,\infty)\end{equation}
in the sense of viscosity solutions.
\hfill $\Box$


\subsection{A lower bound on $Z_*$}

We now derive from equation (\ref{Z_*subsol}) an estimate on $Z_*$.

\begin{lem} \label{lem:compsuper}
One has $Z_* (t,x) \geq \min \{ -t \underline{H}^\star (-x/t),0\}$ for all $(t,x)\in (0,\infty)\times (0,\infty)$, where $\underline{H}^\star$ is the convex conjugate of $\underline{H}$.
\end{lem}

We remind to the reader that the convex conjugate of function $\underline{H}$ is defined by 
$\underline{H}^\star (q) = \sup_{p\in\R} \big( p q - \underline{H} (p)\big)$. 

\bigskip

\noindent {\bf Proof.} 
Define $U(t,x):= -t^{-1} Z_* (t,-tx)$ for all $t,x>0$.  Take $t,x>0$ such that $Z_* (t,-tx) <0,$ then we get from  Proposition \ref{eqZ^*Z_*}:
\begin{equation}\begin{array}{rcl}
\partial_t U(t,x) &=& \displaystyle\frac{1}{t^2} Z_* (t,-tx) -\frac{1}{t} \partial_t Z_* (t,-tx)+\frac{x}{t} \partial_x Z_* (t,-tx) \\
&&\\
&\leq &\displaystyle\frac{-1}{t} U(t,x)  -\frac{1}{t}\underline{H}\big(\partial_x Z_* (t,-tx)\big)+\frac{x}{t} \partial_x Z_* (t,-tx) \\
\end{array}
\end{equation}
in the sense of viscosity solutions. As $\underline{H}(p) + \underline{H}^\star (x) \geq px$ for all $p,x\in\R$, it follows that
\begin{equation}  \label{eqlemcomp}
\partial_t U (t,x)\leq \frac{-1}{t} U(t,x)  + \frac{1}{t}\underline{H}^\star (x).
\end{equation}
On the other hand, the definition (\ref{defsubsuper}) of $Z_*$ yields $Z_* (\alpha t, \alpha x) = \alpha Z_* (t,x)$ for all $(t,x)\in (0,\infty)\times (0,\infty)$  
and $\alpha>0$.
Hence, $U(t,x) = -Z_* (1,-x)$ and in particular, $\partial_t U(t,x) = 0$ in the sense of viscosity solutions for all $(t,x)$ such that $Z_* (t,-tx) <0.$ 
It follows from (\ref{eqlemcomp}) that 
$ U(t,x) \leq \underline{H}^{\star} (x)$. Hence, if $Z_* (t,x)<0$, then:
$$Z_* (t,x) = -t U(t,-x/t) \geq -t \underline{H}^\star (-x/t)\geq \min \{ -t \underline{H}^\star \big(-x/t\big),0\}.$$
If $Z_*(t,x) \geq 0$, then $Z_* (t,x)\geq \min \{ -t \underline{H}^\star \big(-x/t\big),0\}$ is also satisfied.  \hfill $\Box$


\subsection{Conclusion of the proof}

\noindent {\bf Proof of Part 2. of Theorem \ref{spreadingthm}.} Consider $w\in (0,\underline{w})$. Then one has 
$\underline{H}(-p)>(1+\e) pw $ for all $p>0$ and for some $\e>0$. As $\underline{H} (0)>0$ and $\underline{H}$ is continuous from Proposition \ref{prop-Hp}, 
there exists $\delta>0$ such that 
$\underline{H} (-p) >pw+\delta$ for all $p>0$, which means that $-\underline{H}^\star (-w) > 0$. Lemma \ref{lem:compsuper} and the continuity of $\underline{H}^\star$
 yield that for all $x>0$ close to $w$ and $t>0$ close to $1$, one has 
$$Z_* (t,x) \geq \min \{ - t\underline{H}^\star (-x/t),0\}=0.$$
Hence, $(1,w)\in \rm{int} \{Z_* = 0\}$. It follows from Lemma \ref{cvHJ} that
$$\liminf_{\e \to 0} v_\e (1,w) =\liminf_{\e \to 0} u(1/\e, x/\e)=\liminf_{t\to +\infty} u(t,wt) >0.$$
As $\lim_{t\to +\infty} u(t,x)=1$ locally in $x$, 
it follows from Theorems $1.3$ and $1.6$ of \cite{BHNa} that 
$\lim_{t\to +\infty} \sup_{0\leq x \leq w't} |u(t,x)-1|=0$ for all $w'\in (0,w)$. This concludes the proof since $w$ is arbitrarily close to $\underline{w}$. 
 \hfill$\Box$


\section{Application: random stationary ergodic coefficients}

We will need in this Section another notion of principal eigenvalue, introduced in \cite{BerestyckiHamelRossi, BerNirVaradhan}. 
Consider for the moment any (deterministic) second order elliptic operator $\L\phi = a(x) \phi'' + q(x)\phi'+ c(x) \phi$, where $a$, $q$ and $c$ are continuous and uniformly bounded 
functions over $\R$ and $\inf_\R a >0$. 
For all non-empty open interval $I\subset \R$, let
\begin{equation} \label{def-eigenI}
 \Lambda_1 (\L, I) := \inf \{ \lambda , \quad \exists \phi \in\mathcal{C}^2 (I)\cap \mathcal{C}^0 (\overline{I}), \quad \phi >0 \hbox{ in } I, \quad \phi=0 
\hbox{ in } \partial I, \quad \L \phi \leq \lambda \phi \hbox{ in } I\}. 
\end{equation}
This definition seems close to Definition \ref{defgeneigen} except that we do not impose any condition on the limit of $\frac{1}{x} \ln \phi (x)$ as $x\to +\infty$ and 
that we require $\phi = 0$ on the boundary of $I$, which indeed makes a big difference. It is known \cite{BerNirVaradhan} that if $I$ is bounded, then $\Lambda_1 (\L, I)$ is the Dirichlet principal eigenvalue associated with 
$\L$. The properties of $\Lambda_1 (\L, I)$ when $I$ is unbounded have been investigated in \cite{BerestyckiHamelRossi, Rossi2, BerRossipreprint}, where it was proved in 
particular that 
\begin{equation} \label{monLambda} I\subset J \Rightarrow \Lambda_1 (\L, I)\leq \Lambda_1 (\L, J),\end{equation} 
 \begin{equation} \label{limLambda} \Lambda_1 (\L, \R)= \lim_{R\to +\infty} \Lambda_1 \big(\L, B_R (y)\big) \hbox{ for all } y\in\R,\end{equation} 
where $B_R (y) = (y-R,y+R)$. 

\bigskip

Let now turn back to random stationary ergodic coefficients and consider the operators $L_p^\omega$ as in the statement of Theorem \ref{rsethm} for all $\omega \in\O$ 
and $p\in\R$. Theorem \ref{rsethm} will be derived from the following new result, combined with Proposition \ref{prop-eigen} and Lemma \ref{lem-compLambda}. 

\begin{thm}\label{rsethmcorrec}
There exists a measurable set $\O_1$, with $\P (\O_1)=1$, such that for all $p\in\R$ and $\o\in\O_1$:

\begin{enumerate}
 \item one has
$\underline{\lambda_1} (L_p^\omega,\R)=\overline{\lambda_1} (L_p^\omega,\R)$ and this quantity does not depend on $\o\in\O_1$,

\item if $\underline{\lambda_1} (L_p^\omega,\R)> \Lambda_1 (\L^\omega,\R)$, then 
there exists $\phi \in \A_{-\infty}$ such that $L_p^\omega \phi = \underline{\lambda_1} (L_p^\omega,\R) \phi$ in $\R$.
\end{enumerate}
 \end{thm}

Note that part $2.$ of the result is only true for $\underline{\lambda_1} (L_p^\omega,\R)> \Lambda_1 (\L^\omega,\R)$. 
We will prove in Lemma \ref{lem-compLambda} below that $\underline{\lambda_1} (L_p^\omega,\R)\geq \Lambda_1 (\L^\omega,\R)$
for all $p\in\R$ and $\omega\in\O$, but the equality might hold for some $p$ near $0$. 

The existence of classical eigenfunctions (called ``exact correctors'' in the homogenization literature) is an open problem in 
the general framework of nonlinear Hamilton-Jacobi equation (see \cite{LionsSouganidis}). Davini and Siconolfi \cite{DaviniSiconolfi} have proved the existence of 
exact correctors for first order random stationary ergodic Hamilton-Jacobi equations. The framework of Theorem \ref{rsethmcorrec} is different since we consider here a second order 
linear equation, which is a particular second order Hamilton-Jacobi equation. Hence our result gives a second class of equations which admit
exact correctors. Note that Davini and Siconolfi prove the existence of exact correctors as soon as a quantity which plays the role
of an eigenvalue is above a given critical treshold, which corresponds to our constraint 
$\underline{\lambda_1} (L_p^\omega,\R)> \Lambda_1 (\L^\omega,\R)$.

\bigskip

\noindent {\bf Proof of Theorem \ref{rsethmcorrec}.} 1. It has been proved by Nolen \cite{Nolen} when $a\equiv 1$ 
and extended by Zlatos \cite{Zlatos} to general $a$'s that 
there exists a measurable set $\O_1$, with $\P (\O_1)=1$, and a real number 
$\overline{\gamma}>0$ such that for all
$\omega \in\O_1$ and $\gamma>\overline{\gamma}$, there exists a unique positive 
$u=u(\cdot,\omega;\gamma)\in \mathcal{C}^2 (\R)$ which solves
\begin{equation}\label{eq:rseu}
 \big(a(x,\omega) u'\big)' + f_s'(x,\omega,0)u = \gamma u \hbox{ in } \R, \quad u(0,\omega;\gamma)=1, \quad \lim_{x\to +\infty} u(x,\omega;\gamma)=0,
\end{equation} 
while no solution of this equation exists if $\gamma<\overline{\gamma}$. Moreover, for all $\gamma>\overline{\gamma}$ and $\omega \in\O_1$, the limit
$$\mu (\gamma):= \lim_{x\to \pm\infty} \frac{-1}{x} \ln u (x,\omega;\gamma) \quad \hbox{exists and is positive}.$$ 
The function $\mu$ does not depend on $\omega\in\O_1$, it is increasing, concave and converges to $+\infty$ as $\gamma$ tends to $+\infty$. 
Moreover, the quantity $\overline{\gamma}$ does not depend on $\omega$ too and one has 
$\overline{\gamma}= \Lambda_1 (\L^\omega, \R)$ for almost every $\o$. This implies in particular that $\Lambda_1 (\L^\omega,\R)$ is a deterministic quantity. 

\smallskip

2. As $\mu$ is increasing and nonnegative, the limit $\rho_R := \lim_{\gamma\to \overline{\gamma}^+} \mu (\gamma)\geq 0$ exists. 
The function $\mu$ admits an inverse $k: (\rho_R, \infty)\to \big(\overline{\gamma} ,\infty\big)$. 
For all $p>\rho_R$, take $\gamma= k(p)$, consider the solution $u$ of (\ref{eq:rseu}) and let $\phi (x) := u(x) e^{px}$ for all $x \in\R$. This function satisfies
$$L_{-p}^\omega \phi = e^{px} \L^\omega (e^{-px}\phi)=e^{px} \L^\omega u = \gamma e^{px} u = k(p)\phi \hbox{ over } \R.$$

Moreover, one has $\phi>0$ in $\R$ and it follows from the Harnack inequality and elliptic regularity that $\phi'/\phi \in L^\infty (\R)$. 
We know that 
$$\frac{1}{x} \ln \phi (x) = \frac{1}{x} \ln u (x) +p\to -\mu (\gamma) + p= 0 \hbox{ as } x\to \pm\infty.$$
Hence, $\phi \in \A_{-\infty}$ almost surely and it follows from Proposition \ref{prop-eigen} that 
\begin{equation} \label{eq:egvaprse}\forall \omega\in\O_1, \quad\forall p>\rho_R, \quad \underline{\lambda_1} (L_{-p}^\omega ,\R)
=\overline{\lambda_1} (L_{-p}^\omega ,\R)= k(p).\end{equation}
Moreover, as $\rho_R = \lim_{\gamma\to \overline{\gamma}^+} \mu (\gamma)$ and as $p\mapsto \underline{\lambda_1}(L_{-p}^\omega,\R)$ and 
$p\mapsto \overline{\lambda_1}(L_{-p}^\omega,\R)$ are continuous for all $\omega\in\O$, one has 
$$ \forall \omega\in\O_1, \quad \underline{\lambda_1}(L_{-\rho_R}^\omega,\R)=\overline{\lambda_1}(L_{-\rho_R}^\omega,\R)=\overline{\gamma} = \Lambda_1(\L^\omega,\R).$$

Similarly, one can prove the existence of $\rho_L \leq 0$ such that, up to some neglectable modification of $\O_1$, for all $p<\rho_L$ and $\omega \in\O_1$, 
there exists a solution $\phi\in\mathcal{A}_{-\infty}$ and $m(p) > \Lambda_1 (\L^\omega,\R)$ of $L_{-p}^\omega \phi = m(p) \phi$. 
It follows that 
\begin{equation} \label{eq:egvaprse2} \forall \omega\in\O_1, \forall p<\rho_L,\ \underline{\lambda_1} (L_{-p}^\omega ,\R)=\overline{\lambda_1} (L_{-p}^\omega ,\R)= m(p) 
\hbox{ and }  \underline{\lambda_1}(L_{-\rho_L}^\omega,\R)=\overline{\lambda_1}(L_{-\rho_L}^\omega,\R)=\Lambda_1(\L^\omega,\R).\end{equation}

Lastly, the same arguments as in the proof of Proposition 5.1.v of \cite{BerRossipreprint} yield that
$p\mapsto  \overline{\lambda_1} (L_p^\omega ,\R)$ is convex. It follows that  
\begin{equation} \label{eq-ineqlambdaLambda} \forall p\in [\rho_L,\rho_R], \quad \forall \omega\in\O_1, \quad \underline{\lambda_1} (L_{-p}^\omega ,\R)\leq 
\overline{\lambda_1} (L_{-p}^\omega ,\R) \leq \Lambda_1 (\L^\omega,\R).\end{equation}
Hence, $\underline{\lambda_1} (L_{-p}^\omega,\R) > \Lambda_1 (\L^\omega,\R)$ implies $p<\rho_L$ or $p>\rho_R$ and the conclusion follows. 
\hfill $\Box$

\begin{lem} \label{lem-compLambda} There exists a measurable set $\O_2$, with $\P (\O_2)= 1$, such that for all $\omega \in\O_2$ and for all $p\in\R$, one has
 $\underline{\lambda_1} (L_p^\omega, \R) \geq \Lambda_1 (\L^\omega, \R)$.
\end{lem}

\noindent {\bf Proof.} 
Even if it means adding a constant to $f_s'(x,\omega,0)$, we can always assume that 
$\Lambda_1 (\L^\omega,\R)=\overline{\gamma} > 0$ and we thus need to prove that $\underline{\lambda_1} (L_p^\omega,\R)\geq 0$. 

\smallskip

1. For all $R>0$ and $(y,\omega) \in\R\times \O$, let $\Big(\phi_{B_R (y)}^\omega, \Lambda_{1}\big(\L^\omega,B_R(y)\big)\Big)$ the unique eigenelements satisfying
\begin{equation} \label{eq-phiRrse} \left\{ \begin{array}{rclcl}
\L^\omega \phi_{B_R (y)}^\omega &=& \Lambda_{1}(\L^\omega,B_R (y))\phi_{B_R (y)}^\omega  & \hbox{ in }& B_R (y),\\
\phi_{B_R (y)}^\omega &>& 0 &\hbox{ in }& B_R (y),\\
\phi_{B_R (y)}^\omega &=& 0 &\hbox{ over } &\partial B_R (y),\\
\max_{x\in B_R (y)} \phi_{B_R (y)}^\omega (x) &=& 1.&&\\ 
\end{array} \right. \end{equation}
As the eigenelements are continuous with respect to the coefficients $a=a(x,\o)$ and $c=c(x,\omega)$, one can easily check that  
$\omega\in\O\mapsto \Big(\phi_{B_R (y)}^\omega, \Lambda_{1}\big(\L^\omega,B_R(y)\big)\Big)\in \Big(\mathcal{C}^2 \big(B_R (y)\big)\cap \mathcal{C}^0 (\overline{B}_R(y))\Big)\times \R$ are measurable functions for all $y\in\R$ and $R>0$. 

Take $R>0$, $(x,y,z,\omega) \in\R\times\R\times \R\times \O$ and define $\psi (x) := \phi_{B_R (y+z)}^\omega (x+z)$. One has 
$$\begin{array}{rcl} 
\L^{\tau_z\omega} \psi &=& \big( a(x,\tau_z\omega)\psi'(x)\big)'  +f_s'(x,\tau_z \omega,0) \psi (x)= \big( a(x+z,\omega)\psi'(x)\big)' +f_s'(x+z,\omega,0) \psi (x)\\
&&\\
&=& \big( \L^\omega \phi_{B_R (y+z)}^\omega\big) (x+z)= \Lambda_{1}\big(\L^\omega, B_R (y+z)\big)  \phi_{B_R (y+z)}^\omega (x+z) \\
&&\\
&=&\Lambda_{1}\big(\L^\omega, B_R (y+z)\big) \psi (x)\hbox{ in } B_R (y).\\ \end{array} $$
Moreover, one has $\psi>0$ in $B_R (y)$ and $\psi=0$ in $\partial B_R (y)$ and
$$\max_{x\in B_R (y)} \psi (x)=\max_{x\in B_R (y)} \phi_{B_R (y+z)}^\omega (x+z)=1.$$ 
As the solution of (\ref{eq-phiRrse}) is unique, we eventually get $\psi \equiv \phi_{B_R (y)}^{\tau_z \omega} $ and thus 
\begin{equation} \label{eq-rsephiR}
\phi_{B_R (y+z)}^\omega (\cdot+z)\equiv  \phi_{B_R (y)}^{\tau_z \omega} \hbox{ and } \Lambda_{1} \big(\L^\omega,B_R (y+z)\big) =  \Lambda_{1} \big(\L^{\tau_z\omega},B_R (y)\big).\end{equation}
In other words, the eigenelements are random stationary ergodic in $(y,\omega)$. 

\smallskip

2. Next, consider for all $\omega\in\O$ the elliptic equation:
\begin{equation} \label{eq-statu} -\big(a(x,\omega) u'\big)'-2pa(x,\omega)u' = \big(p^2a(x,\omega)+p a'(x,\omega) +f_s'(x,\omega,0) \big)u - u^2  \hbox{ over } \R.\end{equation}
Obviously, $\overline{u}:= \sup_{x\in\R}\big(p^2 a(x,\omega) +p a'(x,\omega)+f_s'(x,\omega,0)\big)$ is a supersolution of this equation. 

On the other hand, for all $(y,\omega)\in\R\times \O$, as $\lim_{R\to +\infty} \Lambda_1 \big(\L^\omega, B_R(y)\big) = \overline{\gamma}$ for all $y\in\R$ 
and $\overline{\gamma}>0$, 
there exists $R(y,\omega)>0$ such that  $\Lambda_1 \big(\L^\omega, B_{R(y,\omega)}(y)\big)= \overline{\gamma}/2$. Moreover, 
as $R\mapsto \Lambda_1 \big(\L^\omega, B_R (y)\big)$ is increasing (see \cite{BerestyckiHamelRossi} for example), $R(y,\omega)$ is uniquely defined and it 
follows from (\ref{eq-rsephiR}) that $R(y+z,\omega)= R(y,\tau_z \omega)$ for all $(y,z,\omega)\in\R\times \R\times \O$. 

It is easy to check that for all $(y,\omega)\in \R\times \O$, the function 
$$\underline{u} (x) = \underline{u} (x,y,\omega):= \left\{ \begin{array}{ccl} \displaystyle
\Lambda_1 \big(\mathcal{L}^\omega, B_{R (y,\omega)(y)}  \big)  \phi_{B_{R(y,\omega)}(y)}^\omega (x)
e^{-p \big(x-y+R(y,\omega)\big)} &\hbox{ if }& x\in B_{R(y,\omega)}(y),\\
&&\\
0 &\hbox{ if }&  x\notin B_{R(y,\omega)}(y),\\
\end{array} \right. $$
 is a subsolution of equation (\ref{eq-statu}). 
Moreover, one has $\underline{u}(x) \leq\Lambda_1 \big(\mathcal{L}^\omega, B_{R (y,\omega)(y)}\big)\leq \overline{\gamma} \leq \overline{u}$. The last inequality is obtained just 
by taking $\phi (x) = e^{px}$ as a test-function in the definition (\ref{def-eigenI}) of $\overline{\gamma}=\Lambda_1 (\L^\omega,\R)$. 
Hence, there exists a minimal solution $u=u(x,\omega)$ of equation (\ref{eq-statu}) in the class of all the solutions satisfying 
$\underline{u} (x,0,\omega)  \leq u (x,\omega)\leq \overline{u}$ for all $(x,\omega)\in\R\times\O$. 

Take $y \in \R$ and let $v(x,\omega):= u(x+y,\tau_{-y}\omega)$ for all $(x,\omega)\in\R\times \O$. The stationarity of the coefficients yields that $v$ is a solution of (\ref{eq-statu}). Moreover, it follows from (\ref{eq-rsephiR}) that $v$ satisfies 
$$\underline{u} (x,0,\omega)=\underline{u} (x+y,0,\tau_{-y}\omega)  \leq u(x+y,\tau_{-y}\omega)=v (x,\omega)\leq \overline{u}$$ 
for all $(x,\omega)\in\R\times\O$. The minimality of $u$ gives $u(x,\omega)\leq v(x,\omega) = u(x+y,\tau_{-y}\omega)$ for all $(x,y,\omega)\in\R\times \R \times \O$. It immediatly follows that 
\begin{equation} \label{eq-rseu} u(x,\tau_{y}\omega) = u(x+y,\omega) \quad \hbox{ for all } \quad (x,y,\omega)\in\R\times \R \times \O.\end{equation}

\smallskip

3. The Harnack inequality and elliptic regularity imply that $x\mapsto \big(u'/u\big) (x,\omega)$ is a bounded function over $\R$ for all $\omega\in\O$. 
As,  $\big(u'/u \big) (x+y,\omega) =\big(u'/u \big) (x,\tau_y\omega)$ for all $(x,y,\omega) \in\R\times \R\times \O$, the Birkhoff ergodic theorem
yields that there exists a measurable set $\O_2\subset \O$ such that $\P (\O_2)=1$ and for all $\omega\in\O_2$:
$$x\mapsto \frac{1}{x}\ln u(x,\omega) \quad \hbox{ converges as } \quad x\to \pm\infty.$$
Moreover, the limits at $+\infty$ and $-\infty$ are equal and do not depend on $\omega \in\O_2$. But as 
$u (x,\omega) \leq \overline{u}$ for all $(x,\omega)$, these limits are necessarily zero: 
$\lim_{|x|\to +\infty} \frac{1}{|x|}\ln u (x,\omega) = 0$ for all $\omega \in\O_2$.
Hence, $u (\cdot,\omega)\in\mathcal{A}_{-\infty}$ for all $\omega \in\O_2$. 
As $L_p^\omega u = u^2 \geq 0$ in $\R$, we can take $u(\cdot,\omega)$ as a test-function in the definition of $\underline{\lambda_1} (L_p^\omega,\R)$, 
leading to $\underline{\lambda_1} (L_p^\omega,\R)\geq 0$ for all $\omega \in\O_2$. 
\hfill $\Box$

\bigskip

\noindent {\bf Proof of Theorem \ref{rsethm}.}
Let $\O_0:= \O_1\cap \O_2$ and define $\rho_R\geq 0\geq \rho_L$, $k(p)$ and $m(p)$ as in Theorem \ref{rsethmcorrec}. 
We know from Theorem \ref{rsethmcorrec} that 
$$\underline{\lambda_1} (L_{-p}^\omega ,\R)=\overline{\lambda_1} (L_{-p}^\omega ,\R) = \left\{ \begin{array}{lcl} m(p) &\hbox{ if }& p<\rho_L\\
                                                                                          k(p) &\hbox{ if }& p>\rho_R\\
                                                                                         \end{array}\right. \hbox{ for all } \omega \in\O_0.$$
Moreover, the proof of Theorem \ref{rsethmcorrec} yields that $\overline{\lambda_1} (L_{-p}^\omega , \R) \leq \Lambda_1 (\L^\omega, \R)$ if $p\in [\rho_L,\rho_R]$. 
Lemma \ref{lem-compLambda} gives $\underline{\lambda_1} (L_{-p},\R)\geq \Lambda_1 (\L^\omega, \R)$ for all $p\in\R$ and thus 
$$\underline{\lambda_1} (L_{-p}^\omega ,\R)=\overline{\lambda_1} (L_{-p}^\omega ,\R) = \Lambda_1 (\L^\omega,\R) \hbox{ if } p\in [\rho_L,\rho_R]  \hbox{ for all } \omega \in\O_0.$$
 \hfill $\Box$

\bigskip

\noindent {\bf Proof of Proposition \ref{rseprop}.} 
Let $\O_0$ as in the statement of Theorem \ref{rsethm}. As 
$$\underline{\lambda_1} (L_p^\omega,\R) = \underline{\lambda_1} \big(L_p^\omega,(R,\infty)\big) = \overline{\lambda_1} \big(L_p^\omega,(R,\infty)\big) 
=\overline{\lambda_1} \big(L_p^\omega,\R\big)$$
for all $\omega \in \O_0$ and $p\in\R$. Hence, it follows from (\ref{defHp-dim1}) and (\ref{defw-dim1}) that $\underline{H}^\omega(p) = \overline{H}^\omega(p)$ and 
$\underline{w}^\omega = \overline{w}^\omega$, which leads to the conclusion. \hfill $\Box$


\section{Application: almost periodic coefficients}
 \label{section-ap}

We prove in this section Theorems \ref{almostperprop} and \ref{almostperthm} and Proposition \ref{asymp-media}.

\bigskip

\noindent {\bf Proof of Theorem \ref{almostperthm}.} The proof relies on a result from Lions and Souganidis \cite{LionsSouganidisap}, 
who proved the existence of approximate correctors in the framework of 
homogenization of Hamilton-Jacobi equations with almost periodic coefficients. Hence, this proof illustrates the strong link between the notion of generalized 
principal eigenvalues used in the present paper and the notion of approximate correctors, which is  used by the homogenization community. 
We will discuss precisely this link, and more generally we will clarify how one can use homogenization techniques to obtain spreading properties, in a forthcoming paper
\cite{BN2}. 

Consider the sequence of equations
\begin{equation}\label{eq-approxap}
a(x) u_\e''+a(x) (u_\e')^2 +q(x) u_\e' +f_s'(x,0)=\e u_\e \hbox{ in } \R.
\end{equation}
As $f_s'(\cdot,0)$ is uniformly bounded, the constants $M_\e = \| f_s'(\cdot,0)\|_{L^\infty (\R)}/\e$ and $-M_\e$ are respectively super and subsolutions 
of (\ref{eq-approxap}). It follows from the Perron's method that there exists a unique solution $u_\e\in\mathcal{C}^{2}(\R)$ of 
equation (\ref{eq-approxap}) such that $-M_\e\leq u_\e \leq M_\e$. In particular, the family $(\e u_\e)_{\e>0}$ is uniformly bounded over $\R$. 

For all $\e>0$, as $u_\e$ is bounded, it follows from (\ref{eq-approxap}) and classical elliptic regularity estimates that $u_\e'$ is bounded. 
Assume that $u_\e'$ reaches a local extremum at $x_0\in\R$. Then (\ref{eq-approxap}) gives
$$a (x_0) \Big( u_\e'(x_0) + \frac{q(x)}{2 a(x)}\Big)^2 = \e u_\e (x_0) - f_s'(x_0,0) + \frac{q(x_0)^2}{4 a(x_0)} \leq \| f_s'(\cdot,0)\|_{L^\infty (\R)}- f_s'(x_0,0) + \frac{q(x_0)^2}{4 a(x_0)} ,$$
from which we easily derive a bound on $u_\e'(x_0)$ which does not depend on $\e$. Distinguishing between the cases where $u_\e'$ reaches its maximum, where it is monotonic 
at infinity and where it is not, we conclude that the family $(u_\e')_{\e>0}$ is uniformly bounded over $\R$. It immediatly follows from (\ref{eq-approxap}) that 
$(u_\e'')_{\e>0}$ is also uniformly bounded over $\R$. 

Next, it has been proved by Lions and Souganidis (see Lemma 3.3 in \cite{LionsSouganidisap}) that, as the coefficients are almost periodic,
$(\e u_\e)_{\e>0}$ converges uniformly in $x\in\R$ as $\e$ tends to $0$. Let $\lambda$ its limit. 
Take $\delta>0$ and $\e>0$ small enough such that $\e u_\e(x)\geq \lambda-\delta$ for all $x\in\R$. Define $\phi:=e^{u_\e}$. 
One has $\phi\in W^{1,\infty}(\R)\cap \mathcal{C}^2(\R)$, $\inf_{\R} \phi>0$ and $\phi'/\phi = u_\e' \in L^\infty (\R)$. Hence, $\phi\in \A_{-\infty}$. 
Moreover, $\phi$ satisfies
$$\L\phi=\e u_\e \phi \geq (\lambda-\delta)\phi \hbox{ in }\R.$$
Hence, one has $\underline{\lambda_1}(\L,\R) \geq \lambda-\delta$ for all $\delta>0$ and thus $\underline{\lambda_1}(\L,\R) \geq \lambda$. 
Similarly, one can prove that $\overline{\lambda_1}(\L,\R\times\R) \leq \lambda$. As $\overline{\lambda_1}(\L,\R\times\R)\geq \underline{\lambda_1}(\L,\R\times\R^N)$, 
this gives the conclusion. \hfill$\Box$

\bigskip

\noindent {\bf Proof of Theorem \ref{almostperprop}.}
Theorem \ref{almostperthm} and (\ref{explicitLp})
give $\overline{\lambda_1}(L_p,\R)=\underline{\lambda_1}(L_p,\R)$ for all $p\in\R$. Thus, using similar arguments as for homogeneous coefficients, one gets
$\underline{H}(p)=\overline{H}(p)=\overline{\lambda_1}(L_p,\R).$
This concludes the proof. \hfill$\Box$

\bigskip

\noindent {\bf Proof of Proposition \ref{asymp-media}.}
As $L_p^*$ is associated with almost periodic coefficients, Theorem \ref{almostperthm} gives $\overline{\lambda_1}(L_p^*,\R)=\underline{\lambda_1}(L_p^*,\R)$.
This implies for all $R>0$:
$$\overline{\lambda_1}(L_p^*,\R)=\overline{\lambda_1}\big(L_p^*,(R,\infty)\big)=\underline{\lambda_1}(L_p^*, (R,\infty)).$$

If only $q$ and $f_s'(\cdot,0)$ were almost periodic at infinity, that is, if $a\equiv a^*$, then Proposition \ref{dep-q} and the convergence (\ref{cv-asympap})
would imply 
$$\overline{\lambda_1}\big(L_p,(R,\infty)\big)-\overline{\lambda_1}\big(L_p^*,(R,\infty)\big)\to 0 \hbox{ as } R\to +\infty$$
and a similar convergence for $\underline{\lambda_1}$. We would then be able to conclude that 
\begin{equation} \label{cv-apinfty}\lim_{R\to+\infty}\overline{\lambda_1}(L_p, (R,\infty))=\lim_{R\to+\infty}\underline{\lambda_1}\big(L_p,(R,\infty)\big)=
\overline{\lambda_1}(L_p^*,\R).\end{equation}

If $a\not \equiv a^*$, then one cannot apply Proposition \ref{dep-q} and an additional argument is needed. As $q^*$ and $c^*$ are arbitrary, we assume that $p=0$ with no 
loss of generality in order to simplify the notations. 
As in the proof of Proposition \ref{almostperthm}, define $(u_\e)_{\e>0}$ by
\begin{equation}\label{eq-approxap2}
a^*(x) u_\e''+a^*(x) (u_\e')^2 +q^*(x) u_\e' +c^*(x)=\e u_\e \hbox{ in } \R.
\end{equation}
We know that there exists a constant $C>1$ such that 
$$\|u_\e'' \|_\infty +\| u_\e'\|_\infty +\|\e u_\e\|_\infty \leq C \hbox{ for all } \e>0,$$
and the proof of Proposition \ref{almostperthm} yields that $(\e u_\e (x))_{\e>0}$ converges to a limit $\lambda$ as $\e\to 0$ uniformly with respect to $x\in\R$,
where $\lambda=\underline{\lambda_1} (\L^*,\R)$. 
Take $\delta >0$ and $\e>0$ small enough so that $|\e u_\e (x) -\lambda | \leq \delta$ for all $x\in\R$ and $R$ large enough so that 
$$\sup_{x>R} \big( |a^*(x)-a(x)| + |q^*(x)-q(x)|+|c^*(x)-f_s'(x,0)|\big) \leq \frac{\delta}{C(1+C)}.$$
Define $\phi_\e := e^{u_\e}$. One has $\phi_\e' = u_\e' \phi_\e$ and $\phi_\e'' = u_\e''\phi_\e + (u_\e')^2 \phi_\e$. 
Hence, this function satisfies
$$a^*(x) \phi_\e'' +q^*(x) \phi_\e' +c^*(x)\phi_\e=\e u_\e \phi_\e \quad \hbox{ in } \R.$$
For all $x>R$, we thus compute
$$\begin{array}{l}
   a(x) \phi_\e'' + q(x) \phi_\e' + f_s'(x,0) \phi_\e \\
\\
 \geq   a^*(x) \phi_\e'' + q^*(x) \phi_\e' + c^*(x) \phi_\e \\
\\
\quad - \big( \sup_{x>R} |a^*(x)-a(x)| \phi''_\e (x)
+\sup_{x>R} |q^*(x)-q(x)| \phi'_\e (x)+\sup_{x>R} |c^*(x)-f_s'(x,0)| \phi_\e (x)\big)\\
\\
\geq  \e u_\e (x) \phi_\e (x)- C (1+C)\sup_{x>R} \big( |a^*(x)-a(x)| + |q^*(x)-q(x)|+|c^*(x)-f_s'(x,0)|\big)\phi_\e (x)\\
\\
 \geq  (\lambda -2\delta) \phi_\e (x).\\
  \end{array}$$
Hence, using $\phi_\e$ as a test-function in the definition of $\underline{\lambda_1}\big(\L,(R,\infty)\big) $, we conclude that for any $\delta>0$, there exists 
some $R$ such that 
$\underline{\lambda_1}\big(\L,(R,\infty)\big) \geq \lambda -2\delta$. Hence, $\lim_{R\to +\infty} \underline{\lambda_1}\big(\L,(R,\infty)\big) \geq \lambda$. 
Similarly, one can prove that $\lim_{R\to +\infty} \overline{\lambda_1}\big(\L,(R,\infty)\big) \leq \lambda$. 
Hence, 
$$\lim_{R\to +\infty} \overline{\lambda_1}\big(\L,(R,\infty)\big) =\lim_{R\to +\infty} \underline{\lambda_1}\big(\L,(R,\infty)\big) = \lambda.$$
Moreover, we know from the proof of Proposition \ref{almostperthm} that $\lambda = \underline{\lambda_1}(\L^*,\R)=\underline{\lambda_1}(\L^*,\R)$. 
As $q^*$ and $c^*$ are abritrary almost periodic functions, we derive from (\ref{explicitLp}):
$$\lim_{R\to +\infty} \overline{\lambda_1}\big(L_p,(R,\infty)\big) =\lim_{R\to +\infty} \underline{\lambda_1}\big(L_p,(R,\infty)\big) = 
\underline{\lambda_1}(L_p^*,\R)=\underline{\lambda_1}(L_p^*,\R)$$
that is, $\overline{H}(p)=\underline{H}(p)=\overline{\lambda_1}(L_p^*,\R)$ for all $p\in\R$. The conclusion immediatly follows. \hfill$\Box$



\end{document}